





\documentstyle{amsppt}
\TagsOnLeft
\overfullrule = 0pt 
\magnification\magstep1

\NoRunningHeads

\hsize=6.0truein
\hoffset=0pt
\voffset=.5truecm
\vsize=8truein

\def\Rn#1{\text{{\it I\kern -0.25emR}$\sp{\,{#1}}$}}
\def\R{\Rn{}}
\def \Lip{\text{\rm Lip\,}}
\def \enk{\{ e_{nk} \}_{k=1}^{\infty}}
\def \ukj{\{ u\ss{k,j} \}_{j=1}^{\infty}}
\def \Uf{f\ss{\Cal U}}
\def \UZ{Z\ss{\Cal U}}

\def \qed {\vrule height6pt  width6pt depth0pt}
\def \e{\epsilon }
\def \d{\delta}
\def\lam{\lambda}
\def\gam{\gamma}
\def\Gam{\Gamma}
\def \nm#1{\left\|#1\right\|}
\def \UX{X\ss{\Cal U}}
\def \UY{Y\ss{\Cal U}}
\def\U{{\Cal U}}
\def \xn{\{ x_n \}_{n=1}^{\infty}}

\def\dstyle{\displaystyle}
\def\ss#1{\lower3pt\hbox{${\scriptstyle #1}$}}

\def\Tt{{\Cal T}^{2}}

\topmatter
\title Affine Approximation of Lipschitz Functions
 and  Nonlinear Quotients\endtitle
 \author S. Bates$^{*+}$, W.B.~Johnson$^{\dag\ddag}$, J. 
Lindenstrauss$^{\ddag+}$, D. Preiss and G. Schechtman$^{\ddag+}$
\endauthor
\thanks  $\!\!\!\!\!\!\!\!^*$Supported by an NSF Postdoctoral Fellowship in 
Mathematics\hfil\break
$^{\dagger}$Supported in part by NSF DMS-9623260\hfil\break
$^{\ddagger}$Supported in part by the
U.S.-Israel Binational Science Foundation\hfil\break
$^+$Participant, Workshop in
Linear Analysis and Probability, Texas A\&M University\endthanks
\address \hfil\break S. Bates, Department of Mathematics, Columbia University,
New York, NY 10027 U.S.A.,\hfil\break
{\tt email: smb\@math.columbia.edu}\hfil\break
W. B. Johnson, Department of Mathematics, Texas A\&M University,
College Station, TX  77843--3368 U.S.A.\hfil\break
{\tt email: johnson\@math.tamu.edu}\hfil\break
J. Lindenstrauss, Institute of Mathematics,
The Hebrew University of Jerusalem, Jerusalem, Israel\hfil\break
{\tt email: joram\@math.huji.ac.il}\hfil\break
D. Preiss, Department of Mathematics,
University College London, London, Great Britain,\hfil\break
{\tt email: dp\@math.ucl.ac.uk}\hfil\break
G. Schechtman, Department of Theoretical Mathematics,
The Weizmann Institute of Science, Rehovot, Israel,\hfil\break
{\tt email: gideon\@wisdom.weizmann.ac.il} \endaddress
\abstract
 New concepts related to approximating a Lipschitz
function between Banach spaces by affine functions are introduced.
Results which clarify when such approximations are possible are
proved and in some cases a complete characterization of the spaces
$X$, $Y$ for which any Lipschitz function from $X$ to $Y$ can be so
approximated is obtained.  This is applied to the study of
Lipschitz and uniform quotient mappings between Banach spaces. It
is proved, in particular, that any Banach space which is a uniform
quotient of $L_p$, $1<p<\infty$, is already isomorphic to a linear
quotient of $L_p$.

\endabstract
\endtopmatter

\document

\footline={\hss\tenrm\number\pageno\hss}

\head 1. Introduction\endhead

In the framework of geometric nonlinear functional analysis there is
by now a quite developed theory of bi-uniform and bi-Lipschitz
homeomorphisms between Banach spaces as well as bi-uniform and
bi-Lipschitz embeddings. The existence  of a bi-uniform or
bi-Lipschitz homeomorphism between $X$ and $Y$ often (although  not
always) implies the existence of a linear isomorphism between these
two spaces (see e.g. [Rib1], [HM], [JLS] and for a complete survey of
the available information the forthcoming book [BL]). The situation
is similar for  bi-Lipschitz embeddings, while for bi-uniform
embeddings the situation is different (see e.g. [AMM] for bi-uniform
embeddings into Hilbert spaces).

The original purpose of the research reported in
this paper was to initiate the study of nonlinear
quotient mappings in the same context. In the
process, we encountered the need to develop new
notions of approximating  Lipschitz functions by
affine ones which go beyond derivatives. We
believe these may prove to be fundamental notions
whose interest goes beyond the particular
applications to nonlinear quotient mappings we
have here. Readers who are interested mainly in
these notions can concentrate on the second half
of this introduction and on section 2.

We consider two related notions of nonlinear quotient mappings.
A uniformly
continuous mapping $F$ from a metric space $X$ onto a metric space
$Y$ is a {\sl uniform quotient mapping} if for each
$\epsilon>0$ there is a $\delta=\delta(\epsilon)>0$ so that for every
$x\in X$,
$F(\text{\rm B}_{\epsilon}(x))\supset \text{\rm B}_{\delta}(Fx)$.  If in
addition the mapping is Lipschitz and $\delta$ can be chosen to be
linear in $\epsilon$   $F$ is called a {\sl Lipschitz quotient
mapping}.  A map
$F: X\to Y$ is a uniform quotient map if and only if $F\times F:
X\times X \to Y \times Y$ maps the uniform neighborhoods of the
diagonal in $X\times X$ onto the set of uniform neighborhoods of the
diagonal in $Y\times Y$.

Linear quotient mappings between Banach spaces are Lipschitz quotient
mappings and bi-Lipschitz (resp. bi-uniform) homeomorphisms are
Lipschitz (resp. uniform) quotient mappings. The class of Lipschitz
or uniform quotients is larger than the class of maps that can be
obtained as compositions of these two obvious classes of examples.
Indeed, for mappings which are composition of maps from the classes
above the inverse image of a point is always connected while, for
example, the map $f(re^{i\theta})=re^{i2\theta}$ from $\Rn{2}$ onto
itself  is a Lipschitz quotient mapping such that the inverse image
of any nonzero point consists of exactly two points. Based on this
simple example one can also build infinite dimensional examples.

As we shall see below (in sections 3 and 4), it is sometimes quite
delicate to check that a given mapping is a Lipschitz or uniform
quotient mapping.

One of the first questions one would like to study about these new
notions of nonlinear quotients is to what extent they can be
``linearized". The simplest question in this direction is when does
the existence of a Lipschitz or uniform quotient mapping from $X$
onto $Y$ imply the existence of a linear quotient mapping.
 
It turns out that existing examples concerning bi-uniform and
bi-Lipschitz homeomorphisms show that linearization is, in general,
impossible for quotient mappings. On the other hand, we show that at
least some of the positive results concerning homeomorphisms can be
carried over to the quotient setting.

In the study of bi-uniform homeomorphisms between Banach spaces, the
first step towards linearization is usually to pass to a bi-Lipschitz
homeomorphisms between ultrapowers of the two Banach spaces. This
step carries over easily to the quotient setting (see  Proposition
3.4 below). We are thus reduced to the question of linearizing
Lipschitz quotient mappings (and of passing back from ultrapowers to
the original spaces). The most natural way to pass from Lipschitz
mappings to linear ones is via differentiation. Recall that a mapping
$f$ defined on open set $G$ in a Banach space
$X$ into a Banach space
$Y$ is called {\sl G\^ateaux differentiable\/}  at $x_0\in G$ if for
every
$u\in X$ $$\lim_{t\to 0} (f(x_0+tu) - f(x_0))/t = D_f(x_0)u$$ exists
and $D_f(x_0)$ (= the differential of $f$ at $x_0$) is a bounded
linear operator from $X$ to $Y$.

The map $f$ is said to be {\sl Fr\'echet differentiable\/} at $x_0$
if the limit above exists uniformly with respect to $u$ in the unit
sphere of $X$.

 Fr\'echet derivatives are very good for linearization purposes when
they are available, but  their existence is rarely ensured. Rather
general existence theorems are known for G\^ateaux derivatives and
these are very useful in the theory of bi-Lipschitz embeddings and
bi-Lipschitz homeomorphisms. Unfortunately, the existence of a
G\^ateaux derivative does not seem to help in the setting of
Lipschitz quotients. As we shall see in Proposition~3.11, the
G\^ateaux derivative of a Lipschitz quotient mapping can  be
identically zero at some points. It is still maybe possible to use
G\^ateaux derivatives in this theory, but Proposition~3.11 shows that
one should take care of choosing carefully a specific point of
differentiability (or maybe a generic one). We do not know how to do
it.
 
As we mentioned above, Fr\'echet derivatives work very nicely when
they exist, which is rare. Even some weaker versions of Fr\'echet
derivatives would suffice for most purposes in our context. One such
version is the 
$\epsilon$-Fr\'echet derivative. A mapping   $f$ defined on open set
$G$ in a Banach space
$X$ into a Banach space
$Y$ is said to be {\sl $\epsilon$-Fr\'echet differentiable} at $x_0$
if there is a bounded linear operator $T$ from $X$ to
$Y$ and a
$\delta>0$ so that
$$\|f(x_0+u) - f(x_0) - Tu\| \le \epsilon\|u\|\quad \text{for}\quad
\|u\| \le
\delta.$$ Unfortunately, there is no known existence theorem even for
points of 
$\epsilon$-Fr\'echet differentiability for Lipschitz mappings in 
situations that are most relevant to us (e.g., $Y$ is infinite
dimensional and isomorphic to a linear quotient of
$X$). There are such theorems if $Y$ is finite dimensional ([LP]) and
also in some situations where $Y$ is infinite dimensional but every
linear operator from $X$ to $Y$ is compact ([JLPS2]). The theorem in
[LP] enables one to prove a result on the ``local" behavior of
uniform quotient spaces which is analogous to a result of Ribe
([Rib1]) on bi-uniform homeomorphisms. Fortunately, it turned out
that in our context the same can be achieved by a linearization
theorem for Lipschitz mappings into finite dimensional spaces whose
proof is considerably simpler than that of the theorem from [LP]. 

Section 2 is devoted to the examination of two new notions of
approximating Lipschitz functions by affine functions. We say that a
pair of Banach spaces $(X,Y)$ has the {\sl approximation by affine
property} (AAP) if for every Lipschitz function $f$ from the unit
ball $B$ of $X$ into $Y$ and every $\epsilon>0$ there is a ball
$B_1\subset B$ of radius $r$, say, and an affine function $L:X\to Y$
so that
$$ ||f(x)-Lx||\le \epsilon r, \quad x\in \text{\rm B}_1.
$$ This inequality is clearly satisfied if $f$ has an
$\epsilon$-Fr\'echet derivative at the center of $B_1$. AAP is
however a definitely weaker requirement than $\epsilon$-Fr\'echet
differentiability. The second property we examine in section 2 is a
uniform version of this property. We say that the pair $(X,Y)$ has
the {\sl uniform approximation by affine property} (UAAP) if the
radius
$r=r(\epsilon,f)$ of
$B_1$ above can be chosen to satisfy $r(\epsilon,f)\ge c(\epsilon)>0$
simultaneously for all functions
$f$ of Lipschitz constant $\le 1$ ($c(\epsilon)$ of course depends
also on $X$ and $Y$). The existence of Fr\'echet derivatives does not
entail any such uniform estimate, and it is easy to see that in some
situations  the affine approximant
$L$ in the definition of UAAP cannot be a Fr\'echet derivative at any
point, even if such derivatives exist almost everywhere.

The main result of section 2, Theorem~2.7,  gives a complete
characterization of the pairs of spaces which have the UAAP: The pair
$(X,Y)$ has the UAAP if and only if  one of the two spaces is finite
dimensional and the other is super-reflexive (that is, has an
equivalent uniformly convex norm). For the study of nonlinear
quotients we only need  part of this characterization, Theorem~2.3,
which states that
$(X,Y)$ has the UAAP if $X$ is super-reflexive and $Y$ is finite
dimensional. Even the AAP under these conditions would suffice for our
application.  

In Proposition~2.8 (respectively, Proposition~2.9) we also
characterize the spaces
$X$ for which $(X,\R)$ (respectively, $(\R,X)$) has the AAP. One
interesting feature of these characterizations is the following: The
pair $(X,\R)$ has the AAP if (and only if) every Lipschitz
$f:X\to\R$ has a point of Fr\'echet differentiability. On the other
hand
$(\R,X)$ may have the AAP also in a situation where there is a
Lipschitz mapping from $\R$ to
$X$   which fails to have,  for some $\epsilon>0$, even a point of
$\epsilon$-Fr\'echet differentiability.

In section 3 we apply Theorem~2.3 to get that if $Y$ is a uniform
quotient of a super-reflexive $X$, then $Y$ is linearly isomorphic to
a linear quotient of an ultrapower of $X$.  This yields, for example,
that a uniform quotient of a Hilbert space is isomorphic to a Hilbert
space. While most of the positive results in section 3 are of
``local" nature, the section ends with a result on the structure of
Lipschitz quotients which is of ``global" nature and is new even for
bi-Lipschitz homeomorphisms (Theorem~3.18): A Lipschitz quotient of
an Asplund space is Asplund.  

Recall that the {\sl Gorelik principle}  [JLS]  says that a
bi-uniform homeomorphism from $X$ onto $Y$ cannot carry the unit ball
in a finite codimensional subspace of $X$ into a ``small"
neighborhood of an infinite codimensional subspace of $Y$. Moreover,
the same is true for the composition of a bi-uniform homeomorphism
with a linear quotient mapping. This raises the natural question of
whether the same holds for general uniform quotient mappings. It
turns out that this is not the case. Section 4 is devoted to several
examples related to this. For example, it follows from
Proposition~4.1 that there is a uniform quotient mapping from
$\ell_2$ onto itself which sends a ball in a hyperplane to zero.

Proposition 4.1  deals with uniform quotient mappings. We do not know
if a Lipschitz quotient mapping can map a finite codimensional
subspace to a point. Questions of this type are of interest also  in
the finite dimensional setting. Lipschitz quotient mappings from
$\Rn{n}$ to itself which have positive Jacobian almost everywhere are
special cases of quasiregular mappings. See [Ric] for a recent book
on this topic.  The topic of quasiregular mappings is quite
developed. One deep theorem of interest to us, due to  Reshetnyak
(see [Ric,~p.~16]), says that the level sets of a quasiregular
mapping are discrete.  We conjecture that there is a result of a
similar nature for Lipschitz quotient mappings from $\Rn{n}$ to
$\Rn{n}$. In Proposition~4.3 we give an elementary proof that the
level sets of a Lipschitz quotient  mapping from 
$\Rn{2}$ to $\Rn{2}$ are discrete (in [JLPS1] it is shown that they
are even finite).

It is also of interest to study uniform quotient mappings from
$\Rn{m}$ to
$\Rn{n}$. In Proposition~4.2 we give an example of a uniform quotient
mapping from $\Rn{3}$ to $\Rn{2}$ which carries a $2$-dimensional
disk to zero.

A forthcoming paper [JLPS1] contains a detailed study of uniform
quotient mappings from $\Rn{2}$ to itself.  It is shown there, for
example,  that such a map may carry an interval to $0$, but that if
the modulus of continuity
$\Omega(\cdot)$ of the uniform quotient map satisfies
$\Omega(t)=o(\sqrt{t})$ as $t\to 0$, then all level sets of the map
are finite.   The results herein and in [JLPS1] indicate that the
subject of nonlinear quotient mappings between Euclidean spaces is a
promising research area for geometric topology and geometric measure
theory.

Unexplained background can be found in [JLS] and the book [LT].

\head 2. Linearizing Lipschitz Mappings\endhead

We begin with a definition:

\proclaim{Definition 2.1} A pair $(X,Y)$ of Banach spaces is said to
have the {\it approximation by affine property\/} (AAP, in short)
provided that for each ball $B$ in $X$,  every Lipschitz mapping
$f:B\to Y$, and each
$\e>0$, there is a ball $B_1\subset B$ and an affine mapping
$g:B_1\to Y$ so that 
$$
\sup_{x\in B_1} ||f(x) - g(x)|| \le \e r \Lip(f),
$$ where $r$ is the radius of $B_1$.  If there is a constant
$c=c(\e) >0$ so that $B_1$ can always be chosen so that its radius is
at least ${c}$ times the radius of $B$, we say that
$(X,Y)$ has the {\it uniform approximation by affine property\/}
(UAAP, in short). 
\endproclaim 

In Theorem~2.7 we show that a pair  of Banach spaces has the UAAP if
and only if one of the spaces is super-reflexive and the other is
finite dimensional. For applications to nonlinear quotients, the most
important fact (Theorem~2.3) is that
  if $X$ is a uniformly smooth Banach space and $Y$ is a finite
dimensional space, then $(X,Y)$ has the UAAP. Notice that the affine
approximant one obtains to the Lipschitz mapping from this result
cannot always be obtained by differentiation even if the domain of
the Lipschitz mapping is the real line (consider the mapping from
$\R$ to $\R$ which is $(-1)^n$ at the integer $n$ and linear on each
interval $[n,n+1]$).

Note that if every Lipschitz mapping from a domain in $X$ into $Y$
has for each $\e>0$ a point of $\e$-Fr\'echet differentiability, then
$(X,Y)$ has the AAP. Consequently, Theorem~11 in [LP] implies that
$(X,Y)$ has the AAP whenever
$X$ is uniformly smooth and $Y$ is finite dimensional. So we could
avoid Theorem~2.3 in the sequel.  On the other hand, the UAAP seems
to be an interesting property in itself and is more tractable than
the AAP.  Moreover, Theorem~2.3 is much easier to prove than
Theorem~11 in [LP].

\proclaim{Proposition 2.2} If $(X,\R)$ has the UAAP, then   $(X,Y)$
has the UAAP for every finite dimensional space $Y$.
\endproclaim

\demo{Proof} We show first that if $\{f_i\}_{i=1}^n$ are real valued
Lipschitz functions with $\Lip(f_i)\break \le 1$ for $1\le i \le n$ on a
ball $B=\text{\rm B}_r(x_0)$ in $X$ and if $\epsilon>0$, then there are a
constant
$\tilde{c}=\tilde{c}(\epsilon,n,X)$, a ball $\widetilde{B}\subset B$
of radius $s\ge \tilde{c} r$ and affine $g_i:X\to \Rn{}$ so that for
$x$ in $\widetilde{B}$, $|f_i(x) - f_i(x)| \le \epsilon s$ for all
$i$. For $n=1$ this is just the definition of UAAP.  The general case
is proved by induction; we just do the case $n=2$ from which the
general case will be clear.  

Let $c(\epsilon)$ be as in Definition 2.1 and find a ball $B_1\subset
B$ of radius $r_1\ge c(\epsilon c(\epsilon)) r$ and an affine $g_1$ so
that for $x$ in $B_1$, $|f_1(x) - g_1(x)| \le \epsilon c(\epsilon)
r_1$.  Applying again Definition 2.1 we find a ball $B_2\subset B_1$
of radius $r_2\ge c(\epsilon) r_1$ and an affine $g_2$ so that for $x$
in $B_2$, $|f_2(x) - g_2(x)| \le \epsilon r_2$.  Clearly on $B_2$ (as
on $B_1$), $|f_1(x) - g_1(x)| \le \epsilon c(\epsilon) r_1 \le
\epsilon r_2$.  

Now suppose that $F$ is a mapping from a ball $B=\text{\rm B}_r(x_0)$ in
$X$  into a normed space $Y$ of dimension $n$ with $\Lip(F)\le 1$. 
Let $\{y_i,y_i^*\}_{i=1}^n$ be an Auerbach basis for
$Y$; that is, $y_i^*(y_j)=\d_{ij}$ and $||y_i||=1=||y_i^*||$.  By what
we proved above there are real valued affine $\{g_i\}_{i=1}^n$ on $X$
and a ball
$\widetilde{B}\subset B$ of radius $s\ge \tilde{c} r$ so that
$|y_i^*F(x) - g_i(x)|\le {{\epsilon s}\over n}$ for $x$ in
$\widetilde{B}$ and $1\le i \le n$.  The affine function $G:X\to Y$
defined by $Gx=\sum_{i=1}^n g_i(x)y_i$ satisfies for $x$ in
$\widetilde{B}$ \ 
$ ||Fx - Gx||
\le
\sum_{i=1}^n |y_i^*(Fx)-g_i(x)|\le \e s$. \hfill\qed
\enddemo

\proclaim{Theorem 2.3} Suppose $X$ is a uniformly smooth Banach
space.   Then   $(X,Y)$ has the UAAP for every finite dimensional
space $Y$.
\endproclaim

\demo{Proof}  By Proposition~2.2, it is enough to check that $(X,\R)$ has
the UAAP.   So let 
$f$ be a mapping from some ball $B=\text{\rm B}_r(x_0)$  in $X$  into
$\R$ with $\Lip(f)\le 1$.

 By rescaling we can assume, without loss of generality, that the
radius of $B$ is one. Now if $|f(x)-f(x_0)|\le 2\e$ for all $x$ in
$B$, then we can well-approximate $f$ on all of $B$ by the function
which is constantly $f(x_0)$, so assume that this is not the case. 
In particular, we assume that $\Lip(f; 1)\ge \e$, where  for $t>0$
$$
 \Lip(f; t) :=\sup_{||x-y||\ge t} {||f(x)-f(y)||\over ||x-y||}.
$$

Since $X$ is uniformly smooth, we can choose $1/2 > \d > 0$ so that
if $||x||=1$ and $||y||\le \d$, then 
$$ ||x+y||=1+x^*(y)+r(y), \tag2.1$$
 with $|r(y)|\le\e ||y||$, and where $x^*$ is the unique norming
functional for $x$; that is, $||x^*||=1=x^*(x)$.  (The remainder
function $r(\cdot)$ of course depends on $x$.)  Define $k=k(\e)$ by 
$$
\left(1+{{\d\e} }\right)^{k-1}\e\le 1<\left(1+{{\d\e}  }\right)^k\e
$$ and observe that there exists  $d$ with $2\ge 2d \ge 4^{-k}$ so
that
$$
\Lip(f;d/2)< \left(1+{{\d\e} }\right)\Lip(f;2d),
\tag2.2$$
 for otherwise iteration would yield
$$ 1\ge \Lip(f;4^{-k})\ge \left(1+{{\d\e}  }\right)^k
\Lip(f;1)\ge \left(1+{{\d\e}  }\right)^k\e >1.
$$ Let $d'$ be the supremum of all $d$, $2\ge 2d \ge 4^{-k}$, which
satisfy
$(2.2)$. By making a translation of the domain of $f$, we can assume,
in view of $(2.2)$, that for some $z$, $-z$ in $B$, where   
$||z||\equiv d\le d'$ with $d'-d$ as small as we like,  we have: 
$$
\Lip(f;d/2)<(1+\d\e){{f(z)-f(-z)}\over{2d}}.
$$

By adding a constant to $f$, we can assume that
$f(z)=dL=-f(-z)$ for some $0<L\le\Lip(f)\le1$.

Let $z^*$ be the norming functional for ${z\over{||z||}}$. Suppose
that $||w||\le \d d$ with $w$ in $B$, and note that the distance from
$w$ to both $z$ and $-z$ is at least $d/2$. We use
$(2.1)$ to estimate
$|f(w)-Lz^*(w)|$.  First, 
$$\align 
 f(w)&=f(w) - f(-z)+f(-z)\le (1+\d\e)L\nm{z+w}-dL\\
&=  (1+\d\e)Ld \nm{{z\over d }+{w\over d}}-dL
\le (1+\d\e)Ld  \left[ 1+z^*\left({w\over d}\right)+ r\left({w\over
d}\right)\right] -dL\\
 &\le (1+\d\e)Ld \left[ 1+{{z^*\left({w}\right)}\over d}+
{{\e||w||}\over d}\right] -dL \le Lz^*(w)+ \e(2+\d+\d\e)\d L d\\
&\le Lz^*(w)+3\e \d d.
\endalign$$

Similarly, starting from
$$ -f(w)=f(z)-f(w)-f(z),
$$ we see that $f(w)\ge Lz^*(w) -3\e \d d$, so that
$$ |f(w)-Lz^*(w)|\le 3\e \d d
$$ whenever $\nm{w}\le \d d$ and $w$ is in $B$.  This completes the
proof  since
$\text{\rm B}\ss{\d d}(0)\cap B$ contains a ball of radius ${\d\over2}d$
and $d$ is bounded below by a constant depending only on
$\e$ and $\d$.
\hfill\qed
\enddemo

It is clear that if $(X,Y)$ has the UAAP (or AAP) then so does any
pair  $(X_1,Y_1)$ with $X$ isomorphic to $X_1$ and $Y$ isomorphic to
$Y_1$.  Thus we could have stated Theorem~2.3 for super-reflexive
$X$ instead of uniformly smooth $X$; we chose the latter because
uniform smoothness arises naturally in the context in which we are
working.  Moreover, the proportional size of the ball on which the
affine approximation is obtained is uniform over all spaces $X$
having a common modulus of smoothness.

Next we show that the isomorphic version of Theorem~2.3 classifies
those spaces $X$ for which $(X,\R)$ has the UAAP.

\proclaim{Proposition 2.4}   $(X,\R)$ has the UAAP if and only if
$X$ is super-reflexive.
\endproclaim

\demo{Proof} The ``if" direction follows from Theorem~2.3.  For the other
direction, suppose for contradiction that $X$, and hence also
$X^*$, is not super-reflexive.  Given any $\e>0$ and natural number
$N$ we get from [JS] a sequence $\{x_i^*\}_{i=1}^{2^N}$ of unit
vectors in $X^*$ so that for every $m=1,2,...,2^N$ and $a_i\ge 0$,
$$ ||\sum_{i=1}^m a_ix_i^*-\sum_{i=m+1}^{2^N} a_ix_i^*||\ge
\left(1-{\e\over 2}\right)\sum_{i=1}^{2^N} a_i. \tag2.3$$
 Define a dyadic tree as follows:  For $k=0,1,...,N$ and
$j=0,1,...,2^k-1$, set
$$ x_{j,k}^*=2^{k-N} \sum_{i=j\cdot 2^{N-k}+1}^{(j+1)2^{N-k}} x_i^*.
\tag2.4$$
 By $(2.3)$, $\left(1-\e/2\right)\le ||x_{j,k}^*||$,  by
${(2.4)}$,  
$$ x_{j,k}^*={{x_{2j,k+1}^*+x_{2j+1,k+1}^*}\over 2} \tag2.5
$$ 
and by $(2.3)$,
$$ ||x_{2j,k+1}^*-x_{2j+1,k+1}^*||>2(1-\e). \tag2.6
$$

Define $0 < \lam < 1$ by  $\lam^N=1/3$ and define $f$ on $X$ by 
$$ f(x)=\max_{j,k} \lam^k\left(||x||+2|\langle
x_{j,k}^*,x\rangle|\right).\tag2.7
$$
 So $f$ is an equivalent norm on $X$ and $\Lip(f)=3$.  

Suppose $r$ satisfies 
$$
\e r > 18(1-\lam) \quad \text{and} \quad  \text{\rm B}_r(x)\subset \text{\rm
B}_1(0). \tag2.8
$$
 We show that the restriction of $f$ to $\text{\rm B}_r(x)$ is not close
to an affine mapping.  Choose $j$ and $k$ with $k$ as  small as
possible so that
$$ f(x)=\lam^k\left(||x||+2|\langle
x_{j,k}^*,x\rangle|\right),\tag2.9
$$ 
and assume for definiteness that $\langle x_{j,k}^*,x\rangle \ge
0$.  Since $\lam^N=1/3$ and $k$ is minimal, we see that $k<N$.   By
$(2.7)$ and $(2.9)$, 
$$
\lam^{k+1}\left(||x||+2\left[|\langle x_{2j,k+1}^*,x\rangle|\vee
|\langle x_{2j+1,k+1}^*,x\rangle |\right]\right)\le f(x)
$$ so that 
$$\displaystyle |\langle x_{2j,k+1}^*,x\rangle|\vee |\langle
x_{2j+1,k+1}^*,x\rangle| \le
\lam^{-1}\left(\langle
x_{j,k}^*,x\rangle+2^{-1}(1-\lam)||x||\right).$$ Therefore, for
$i=0,1$ we have
$$
\langle x_{2j+i,k+1}^*,x\rangle \ge  (2-\lam^{-1}) \langle
x_{j,k}^*,x\rangle -(2\lam)^{-1}(1-\lam)||x||.
\tag2.10
$$ 
By $(2.3)$ there is $y$ in  $Y$ with $||y||=1$ and 
$\langle {x_{2j,k+1}^*-x_{2j+1,k+1}^*},y\rangle >2(1-\e)$, and hence
$$
\langle x_{2j,k+1}^*, y\rangle > 1-2\e \quad \text{and} \quad 
\langle x_{2j+1,k+1}^*, y\rangle < -1+2\e. \tag2.11
$$
 Therefore, by using $(2.10)$ and $(2.11)$ in the second inequality
below, we get
$$\align
f(x+ry) &\ge  \lam^{k+1}\left(||x+ry||+2\langle
x_{2j,k+1}^*,x+ry \rangle\right)\\
&\ge \lam^{k+1}\left(||x||-r
+2(2-\lam^{-1})\langle x_{j,k}^*,x\rangle -\lam^{-1}  (1-\lam) ||x||
+ 2(1-2\e)r\right)\\
&\ge \lam f(x) + \lam^{k+1}(1-4\e)r-3\lam^k(1-\lam) 
\quad\quad\quad \quad\quad\text{by}\ \ \  (2.9)\\
& \ge f(x) +
3^{-1}(1-4\e)r-6(1-\lam)
\quad\text{since}\ \ \  |f(x)|\le 3 \text{ and }
\lam^{k+1}\ge 3^{-1}\\
& \ge f(x) + 3^{-1}(1-5\e)r.
\quad\quad\quad\quad\quad\quad
\quad\quad\text{by}\
\ \  (2.8)
\endalign$$
 Similarly, $f(x-ry) \ge f(x) + 3^{-1}(1-5\e)r$, and hence
$$
\sup_{y\in\text{\rm B}_r(x)} |f(y)-L(y)|\ge  6^{-1}(1-5\e)r
$$  for every affine function $L$. \hfill\qed
\enddemo

If we reverse the order of $\R$ and $X$ we  get the same
characterization of the  UAAP for the pair of spaces.

\proclaim{Proposition 2.5} A Banach space $X$ is super-reflexive if
and only if 
$(\R,X)$ has the UAAP.
\endproclaim

\demo{Proof} 
Assume first that $X$ is super-reflexive. Without loss of
generality we may assume it is actually uniformly convex. In
Proposition~2.6 we give a soft proof that $(\R,X)$ has the UAAP, but
here we give a proof in the spirit of the proof of Theorem~2.3 which
can provide a specific estimate of the size of the ball upon which
the affine approximation is valid in terms of the modulus of
convexity of $X$ and the degree of approximation  (however, we do not
actually make the estimate).

Let
$f:I\to X$ be a function with
$\Lip(f)\le 1$ defined on an interval $I$ which, without loss of
generality, we assume is of length $1$.  Given $\e$, we have a largest
positive $d$ with   $d\le 1$  
$$
\Lip(f;d/4)\le (1+\e/2)\Lip(f;d).\leqno\text{\rm and}
$$ 
As in the proof of Theorem~2.3, $d$ is  bounded away from zero by a
function of $\e$.  By translations in the domain and range, we may
assume that
$0$, $d\in I$, 
$$
\Lip(f;d/4)<(1+\e){{\|f(d)-f(0)\|}\over{d}},
$$
$f(0)=0$, and $f(d)=dx$  for some $x\in X$ with
$0<\|x\|\le\Lip(f)\le1$.

If $d/4\le r<3d/4$,
$$
\|f(r)- dx\|=\|f(r)- f(d)\|\le (d- r)\Lip(f;d/4)\le (d- r)(1+\e)\|x\|
$$
 $$
\|f(r)\|=\|f(r)- f(0)\|\le  r(1+\e)\|x\|. \leqno{\text{and}}
$$

It follows that
$$
\|f(r)\|+\|dx-f(r)\|\le
(1+\e)d\|x\|=(1+\e)\|f(r)+dx-f(r)\|.\tag2.12
$$

We shall denote below by $\d(\e)$ a positive function, not the same in
every  instance, depending only on the modulus of uniform convexity
of $X$, and  which tends to zero as $\e$ tends to zero. A simple
consequence of the uniform convexity definition is that if
$\|x\|+\|y\|\le(1+\e)\|x+y\|$  for two vectors  in $X$ for which
${1\over 10}
\le {{\|x\|}\over{\|y\|}}\le 10$ then 
$\left\|{x\over{\|x\|}}-{y\over{\|y\|}}\right\|\le\d(\e)$. When
applied to (2.12),  this implies  that $\|f(r)-ax\|\le \d(\e)d\|x\|$
for some scalar
$a$ depending on $r$. Considering now the real valued Lipschitz
function
$\|f(r)\|$, we easily get that, for $d/4\le r\le 3d/4$,
$\big|\|f(r)\|-r\|x\|\big|\le \d(\e)d\|x\|$. Combining inequalities
we get that
 $\|f(r)-rx\|\le \d(\e)d\|x\|\le \d(\e)d $ for $d/4\le r\le 3d/4$.

If $X$ is not super-reflexive we can assume, by [JS], that $X$
contains, for each positive integer $n$ and for each $\e>0$, a
normalized sequence 
$\{e_i\}_{i=1}^n$ satisfying \hfill\break
\text{$\|\sum_{i=1}^k a_ie_i-\sum_{i=k+1}^n a_ie_i\|
\ge(1-\e)\sum_{i=1}^n a_i$,}  for all $k=1,\dots,n$ and all $a_i\ge0$,
$i=1,\dots,n$. Fix $n$ and define 
$f:[0,n]\to X$ by
$$ f(t)=\sum_{i=1}^{[t]}e_i +(t-[t])e_{[t]+1}.
$$ Then f is Lipschitz with constant $1$ and $f^{-1}$ is Lipschitz
with constant $(1-\e)^{-1}$. Assume $k$ and $l$ are positive integers
with  
$k+2l\le n$ and $A:[k,k+2l]\to X$ is an affine map which approximates
$f$ to within $l/4$ on $[k,k+l]$. Then, $\|A(k)-\sum_{i=1}^k e_i\|\le
l/4$,
\text{
$\|A(k+2l)-\sum_{i=1}^{k+2l}e_i\|\le l/4$},  and 
$\|A(k+l)-\sum_{i=1}^{k+l}e_i\|\le l/4$. It then follows that
$$ (1-\e)l\le{1\over 2}
\|\sum_{i=k+1}^{k+l}e_i-\sum_{i=k+l+1}^{k+2l}e_i\|=\|\sum_{i=1}^{k+l}e_i-{1\over
2}(\sum_{i=1}^k e_i+\sum_{i=1}^{k+2l}e_i)\|\le l/2
$$ which, for $\e<1/2$, contradicts the assumption on
$\{e_i\}_{i=1}^n$.
\hfill\qed
\enddemo

\proclaim{Proposition 2.6} Assume that  $X$ is super-reflexive and
$Y$ is finite dimensional. Then $(Y,X)$ has the UAAP.
\endproclaim

\demo{Proof}
 Suppose that $(Y,X)$ fails  the UAAP and let $\U$ be a free
ultrafilter on the natural numbers. We claim  that then
$(Y,X_\U)$ fails the AAP, where $X_\U$ denotes the ultrapower of $X$
with respect to $\U$.   Having established the claim, we complete the
proof by pointing out that   $X_\U$ is reflexive since $X$ is
super-reflexive, and hence every Lipschitz mapping from a ball in $Y$
to
$X_\U$ has a point of differentiability; in particular, $(Y,X_\U)$ has
the AAP.
\enddemo

We now prove the claim.   Since $(Y,X)$ fails  the UAAP, it is easy to
see that there exist
$\e>0$ and mappings $f_n$ from $\text{\rm B}_1^Y(0)$ into $X$ with
$\Lip(f_n)=1$, $f_n(0)=0$, so that for all  balls $B\subset \text{\rm
B}_1^Y(0)$ of radius $r=r(B)$ at least $1/n$ we have for all affine
mappings $L:Y\to X$ the estimate
$$ ||f_n-L||_B\ge \e r, \tag2.13
$$ 
where $\displaystyle ||g||_C:= \sup_{y\in C} ||g(y)||$. Let $f_\U$
be the ultraproduct of the mappings $f_n$, defined for ${y}$ in $\text{\rm
B}_1^{Y}(0)$ by $f_\U({y})=(f_n(y))$.  Let $B\subset \text{\rm B}_1^Y(0)$
be a ball and let $r$ be the radius of $B$.  Now if $L:Y\to X_\U$ is
affine, then since $Y$ is finite dimensional, $L$ is the ultraproduct
of some sequence $(L_n)$ of affine mappings from $Y$ to $X$.  But
then by
$(2.13)$, for all balls $B\subset \text{B}_1^Y(0)$, $\displaystyle
||f_\U-L||_B=\lim_{n\in U} ||f_n-L_n||_B
\ge \e r$.  This means that $(Y,X_\U)$ fails the AAP.\hfill\qed

We now  characterize those pairs of Banach spaces which have the UAAP.

\proclaim{Theorem 2.7} Let $X$ and $Y$ be nonzero Banach spaces. The
pair
$(X,Y)$  has the UAAP if and only if one of the spaces is
super-reflexive and the other is finite dimensional.
\endproclaim

\demo{Proof} 
 If one of the spaces is super-reflexive and the other is
finite dimensional, then $(X,Y)$  has the UAAP by Theorem~2.3 or by
Proposition~2.6.  So assume that the spaces are nonzero and  $(X,Y)$ 
has the UAAP.  It is easy to see that if $X_0$ is a complemented
subspace of $X$ and $Y_0$ is a complemented subspace of $Y$, then
$(X_0,Y_0)$ has the UAAP.  Therefore,
$(X,\R)$ and
$(\R,Y)$ have the UAAP and hence, by Proposition~2.4 and
Proposition~2.5, 
$X$ and $Y$ are super-reflexive.  Now if both are infinite
dimensional, then   both contain uniformly complemented copies of
$\ell_2^n$ for all $n$ ([FT], [Pis]).  Although
$(\ell_2^n,\ell_2^n)$ has the UAAP for each fixed $n$, the estimate
for the size of the ball on which a Lipschitz constant one mapping
from
$\text{B}_1^{\ell_2^n}(0)$ into
$\ell_2^n$ must have a linear approximate within a given error tends
to zero as $n\to \infty$. Indeed,  consider the mapping which takes
$\{a_i\}_{i=1}^n$ to
$\{|a_i|\}_{i=1}^n$. Since any $x\in \text{B}_1^{\ell_2^n}(0)$ has a
coordinate whose absolute value is $\le 1/\sqrt{n}$, any ball of
radius $r>2/\sqrt{n}$ contained in $\text{B}_1^{\ell_2^n}(0)$ contains
a segment of the form $x+(-r/2,r/2)e_i$ for some basis vector $e_i$.
Clearly, on this segment, the mapping above cannot be approximated by
an affine mapping to a degree better than $r/4$.   This easily
implies that 
$(X,Y)$  fails the UAAP if both spaces contain uniformly complemented
copies of $\ell_2^n$ for all $n$.
\hfill\qed
\enddemo

For  examples of pairs of infinite dimensional spaces for which 
$(X,Y)$ has the AAP  we refer to [JLPS2].  In the
rest of this section we restrict attention to the
situation in which either $X$ or $Y$ is the scalar
field.

\proclaim{Proposition 2.8} The following are equivalent for a Banach
space
$X$.
\itemitem{(i)} $(X,\R)$ has the AAP.
\itemitem{(ii)} $X$ is Asplund.
\itemitem{(iii)} Every real valued Lipschitz mapping from a domain in
$X$ has a point of Fr\'echet differentiability.
\endproclaim

\demo{Proof}  Recall that $X$ is  Asplund provided every real valued convex
continuous function on a convex domain $U$ in $X$ is Fr\'echet
differentiable on a dense $G_\d$ subset of $U$.  Various equivalents
to this  can be found e.g. in section I.5 of [DGZ].

That condition (iii) implies both (i) and (ii) is clear. In [Pre] it
is proved that (ii) implies (iii). The  implication (i) $\Rightarrow
$  (ii) is a simple consequence of a result of Leach and Whitfield
[LW] (or see  Theorem~5.3 in [DGZ]).  They proved  that if $X$ is not
Asplund,  then there is an equivalent norm (which we might as well
take to be the original norm 
$||\cdot||$) on $X$ and a 
$\d>0$ so that every weak$^*$-slice of the dual ball in $X^*$ has
diameter larger than $\d$.  This  implies that for every $x$ in
$X$ and
$\e>0$  the diameter of the set $\{x^*\in \text{B}_1^{X^*}(0) :  
x^*(x)> ||x||-\e
\}$ is larger than $\d$.  Their further  computations used to prove
that
$||\cdot ||$ is a so-called rough norm on $X$ also yield that
$(X,\R)$ fails the AAP.  Explicitly, it is enough to check that if
$x$ is in
$X$  and $r>0$, then   there is a vector $u$ in $X$ with $||u||=r$ so
that
$||x+u||+||x-u||-2||x|| > {\d\over 2}r$.  To do this,  choose norm one
linear functionals $x_1^*$,
$x_2^*$ so that
$x_i^*(x)>||x||-{{\d r}\over 4}$ for $i=1,2$ but
$||x_1^*-x_2^*||>\d$. So we can take $u$ in $X$ with $||u||=r$ and
$(x_1^*-x_2^*)(u) > \d r$.  Putting things together, we see that
$$ ||x+u||+||x-u|| \ge  x_1^*(x+u)+x_2^*(x-u) > 2||x|| - {{\d r}\over
2} +
\d r,
$$ as desired.\hfill\qed
\enddemo

 It is also possible to characterize the spaces $X$ for which
$(\R,X)$ has the AAP. Let $\{{\Cal F}_n\}_{n=0}^\infty$ be a sequence
of $\sigma$-fields on $[0,1)$ with ${\Cal F}_0=\{\emptyset,[0,1)\}$,
${\Cal F}_n$ is generated by exactly
$2^n$ intervals of the form $[a,b)$, and the atoms of ${\Cal
F}_{n+1}$ are  obtained from those of ${\Cal F}_n$ by splitting each
of the intervals generating 
${\Cal F}_n$ into two subintervals.

An $X$ valued martingale $\{M_n\}$ with respect to such a sequence of
$\sigma$-fields is called a {\sl generalized dyadic martingale}. If
each atom of ${\Cal F}_n$ has measure exactly $2^{-n}$ we call the
martingale {\sl dyadic}. $\{M_n\}$ is said to be {\sl $\d$-separated}
if $\|M_n-M_{n+1}\|\ge
\d$ a.e. for all $n$. It is {\sl bounded} if
$\sup\|M_n\|_\infty<\infty$. See [KR] or the book [Bou] for
connections between the geometry of $X$ and the existence of
$\d$-separated, bounded, dyadic or generalized dyadic martingales
with values in $X$. It seems not to be known whether the existence of
an
$X$ valued, bounded, 
$\d$-separated generalized dyadic martingale implies the existence of
a dyadic martingale with the same properties (possibly with a
different $\d>0$).

\proclaim{Proposition 2.9} For a Banach space $X$, $(\R,X)$ has the
AAP if and only if for every $\d>0$ there exists no $X$ valued,
bounded, 
$\d$-separated generalized dyadic martingale.
\endproclaim

\demo{Proof} If $\{M_n\}$ is an $X$ valued,  
$\d$-separated generalized dyadic martingale satisfying
$\|M_n\|_\infty\le 1$,  say, define, for each $n$, $f_n:[0,1)\to X$
by $f_n(t)=\int_0^t M_n(s)ds$. Then $\Lip(f_n)\le 1$ for all $n$ and
the sequence $\{f_n(t)\}$ is eventually  constant for each $t$ which
is an end point of one of the intervals generating  one of the ${\Cal
F}_n$'s. The separation and boundedness conditions imply that the set
of such $t$'s is dense in $[0,1]$. Thus $f_n(t)$ converges for all
$t$ in $[0,1)$ to a function $f$ satisfying
$\Lip(f)\le 1$.

If $a<b<c$ with $[a,b)$, $[b,c)$ atoms of ${\Cal F}_{n+1}$ and
$[a,c)$ an atom of
${\Cal F}_n$, then $f(x)=f_{n+1}(x)$ for $x=a,b,c$. Put
$A=M_{n+1}^{[a,b)}$, 
$C=M_{n+1}^{[b,c)}$ and $B=M_{n}^{[a,c)}$, where $M_k^S$ denotes the
constant value of $M_k$ on the atom $S$ of ${\Cal F}_k$. Note first
that
${{c-a}\over{b-a}}\le {2\over\d}$. Indeed,
$B={{b-a}\over{c-a}}A+{{c-b}\over{c-a}}C$ and, if
${{b-a}\over{c-a}}<{\d\over 2}$, then
$\left|1-{{c-b}\over{c-a}}\right|={{b-a}\over{c-a}}<{\d\over 2}$ and
$$
\|B-C\|<{\d\over 2}\|A\|+{\d\over 2}\|C\|<\d,
$$  a contradiction.  

Now if $L$ is an affine map on $[a,c)$ with $L(a)=f(a)$ and
$\|L(x)-f(x)\|\break <(\d^2/4)(c-a)$ on $[a,c)$, then 
$$\align 
\|A-B\|=&\left\|{{f(b)-f(a)}\over {b-a}}-
{{f(c)-f(a)}\over {c-a}}\right\|\\
\le& 
\left\|{{L(b)-L(a)}\over {b-a}}- {{L(c)-L(a)}\over
{c-a}}\right\|+{\d^2\over 4}{{c-a}\over{b-a}} +{\d^2\over 4}\\
=&{\d^2\over 4}\left({{c-a}\over{b-a}}+1\right)\le {\d^2\over
4}\left(1+{2\over\d}\right)<\d,
\endalign$$
  a contradiction. Since every  interval contains an interval of
comparable length which is an atom of one of the ${\Cal F}_n$'s, it
follows  that  $(\R,X)$ does not have the AAP.

If $(\R,X)$ does not have the AAP, let $f:[0,1)\to X$ be a function
with
$\Lip(f)=1$ which for some $\d>0$ has the property that for every
subinterval $[a,c)$ of
$[0,1)$, there exists a  $b$ with $a<b<c$ and 
$$
\left\|f(b)-{{c-b}\over{c-a}}f(a)-{{b-a}\over{c-a}}f(c)\right\|>\d(c-a).
\tag2.14
$$
Then define a generalized dyadic martingale as follows:
$M_0=f(1)-f(0)$ everywhere on $[0,1)$. Assume $M_n$ was already
defined and, for each atom interval $[a,c)$ in ${\Cal F}_n$,
\text{$M_n={{f(c)-f(a)}\over{c-a}}$} on
$[a,c)$. Then, for each such interval pick a $b$ as in $(2.14)$, add
the two intervals $[a,b),[b,c)$ to  ${\Cal F}_{n+1}$, and define
$M_{n+1}={{f(b)-f(a)}\over{b-a}}$ on
$[a,b)$ and $M_{n+1}={{f(c)-f(b)}\over{c-b}}$ on
$[b,c)$. Then 
$$
M_{n}^{[a,c)}={{b-a}\over{c-a}}M_{n+1}^{[a,b)}+{{c-b}\over{c-a}}M_{n+1}^{[b,c)},
$$ and
$$\align
\|M_{n+1}^{[a,b)}-M_{n}^{[a,c)}\|=&
\left\|{{f(b)-f(a)}\over{b-a}}-{{f(c)-f(a)}\over{c-a}}\right\|\\
=&{1\over
{b-a}}\left\|f(b)-(1-{{b-a}\over{c-a}})f(a)-{{b-a}\over{c-a}}f(c)\right\|\\
=&{1\over
{b-a}}\left\|f(b)-{{c-b}\over{c-a}}f(a)-{{b-a}\over{c-a}}f(c)
\right\| >{{\d(c-a)}\over{b-a}}>\d.
\endalign$$
 Similarly, also $\|M_{n+1}^{[b,c)}-M_{n}^{[a,c)}\|>\d$. So
$\{M_n\}$ is an 
$X$ valued,  
$\d$-separated generalized dyadic martingale with $\|M_n\|_\infty\le
1$.
\hfill\qed
\enddemo

 It is well known  that
$X$  fails the Radon-Nikodym property if and only if there is a
bounded,
$\d$-separated, $X$ valued martingale (see e.g. section V.3 in
[DU]).   There is an ingenious example in [BR] of a subspace
$X$ of
$L_1$ which does not have the Radon-Nikodym property, yet for all
$\d>0$ there is no bounded,
$\d$-separated, $X$ valued generalized dyadic martingale. Thus
$(\R,X)$ can have the AAP when $X$ fails the Radon-Nikodym property. 
The characterization of the Radon-Nikodym property in terms of
differentiablility properties is given by the following (mostly known)
proposition:

\proclaim{Proposition 2.10} The following are equivalent for a Banach
space $X$.  
\itemitem{(i)} $X$ has the Radon-Nikodym property.
\itemitem{(ii)} Every Lipschitz mapping from $\R$ into $X$ is 
differentiable almost everywhere.
\itemitem{(iii)}  Every Lipschitz mapping from $\R$ into $X$ has for every
$\e>0$ a point of $\e$-Fr\'echet differentiability.
\endproclaim

\demo{Proof} $(i)\Rightarrow (ii)$ is classical, so one only needs to prove
$(iii)\Rightarrow (i)$. The proof is basically the same as the proof
of the first half of Proposition~2.9: If $X$ fails the Radon-Nikodym
property then there is a bounded,
$\d$-separated, $X$ valued martingale. That is, for some sequence
$\{{\Cal F}_n\}_{n=0}^\infty$ of
$\sigma$-fields on $[0,1)$ 
 such that
 ${\Cal F}_0=\{\emptyset,[0,1)\}$ and each ${\Cal F}_n$ is generated
by intervals, there is an $X$ valued martingale $\{M_n\}$ with respect
to this sequence of  $\sigma$-fields such that $||M_n||_\infty\le 1$
and $||M_n(x)-M_{n+1}(x)||\ge \d$ for all $n$ and all $x\in [0,1)$.
Define $f$ as in the beginning of the proof of Proposition~2.9. If
$x_0$ is a point such that
$$ ||f(x_0+u) - f(x_0) - Lu|| \le \epsilon|u|\quad \text{for}\quad
|u| \le
\eta
$$ for some linear $L$, let $a\le b<c\le d$ with $[b,c)$ an atom of
${\Cal F}_{n+1}$ and $[a,d)$ an atom of
${\Cal F}_n$ with $x_0\in[b,c)$ and $d-a<\eta$. Then
$f(x)=f_{n+1}(x)$ for $x=a,b,c,d$, and 
$$\align
 ||f(c)&-f(b)-L(c-b)||\cr &=||f(x_0+c-x_0)-f(x_0)-L(c-x_0)
-(f(x_0+b-x_0)\\
&-f(x_0)-L(b-x_0))||\\
&\le 2\epsilon(c-b).
\endalign$$ 
Similarly,
$$ ||f(d)-f(a)-L(d-a)||\le 2\epsilon (d-a).
$$ Consequently,
$$ ||M_n(x_0)-M_{n+1}(x_0)||=\left\|{{f(d)-f(a)}\over{ d-a}}-
{{f(c)-f(b)}\over{ c-b}}\right\|\le 4\epsilon,
$$ which is impossible if $\epsilon<\d/4$. \hfill\qed
\enddemo

\noindent{\bf Remark 2.11.} The implication $(ii)\Rightarrow (i)$ in
Proposition~2.10 has been known to experts for a long time.  It is
mentioned without proof in [Ben].  The only published proof of which
we are aware appears in [Qia], where a more general theorem is proved.

\head 3. Nonlinear quotient mappings\endhead

We begin with two definitions. The definition of co-Lipschitz 
appears in Chapter 1, section 1.25 of [Gro].  A mapping $f$ is
co-uniform if and only if the sequence $(f,f,f,\dots)$ is uniformly
open in the sense used in Chapter 10, section 2 of [Why].  Co-uniform
mappings and uniform quotients in the context of general uniform
spaces are discussed in [Jam], where they are used to develop a
theory of uniform transformation groups and uniform covering spaces.

\proclaim{Definition 3.1} A mapping $T$ from a metric space $X$ to a
metric space $Y$ is said to be {\sl co-uniformly continuous\/}
provided that for each $\e>0$ there exists
$\d=\d(\e)>0$ so for every $x$ in $X$, $T\text{\rm B}_{\e}(x)\supset \text{\rm 
B}_{\d}(Tx)$.  If  $\d(\e)$ can be chosen larger than 
 $\e/ C$ for some $C>0$, then $T$ is said to be {\sl co-Lipschitz},
and  the smallest such $C$ is denoted by $\text{co-Lip}(T)$.
\endproclaim

\proclaim{Definition 3.2} A mapping $T$ from a metric space $X$ to a
metric space $Y$ is said to be a {\rm uniform quotient mapping\/}
provided $T$ is uniformly continuous and co-uniformly continuous. If
$T$ is Lipschitz and co-Lipschitz, then $T$ is called a {\rm
Lipschitz quotient mapping}.
\endproclaim

A space $Y$ is said to be a {\sl uniform quotient\/} (respectively, 
{\sl Lipschitz quotient\/}) of a space $X$ provided there is a
uniform quotient mapping (respectively, a Lipschitz quotient mapping)
from $X$ onto $Y$.

The linear theory is simplified by the fact that a surjective bounded
linear operator  between Banach spaces is automatically a quotient
mapping.  Of course, a surjective     Lipschitz mapping from $\R$ to
$\R$ which is a homeomorphism need not be even a uniform quotient
mapping. Moreover, a surjective Lipschitz mapping does not carry any
structure:  in [Bat] it is shown that if $X$ and $Y$ are Banach
spaces with $X$ infinite dimensional, and the density character of
$X$ is at least as large as that of $Y$, then there is a Lipschitz
mapping of $X$ onto $Y$.  In this section we show that, in contrast
to this, uniform and Lipschitz quotient mappings do preserve some
structure.  

A connection between uniformly continuous mappings and Lipschitz
mappings
 is provided by the well-known fact that a uniformly continuous
mapping from a convex domain is ``Lipschitz for large distances". 
There is a ``co" version of this:

\smallskip

\noindent{\bf Remark 3.3.}
 If $T$ is a  uniformly continuous mapping
from a convex set  $X$  onto a  set $Y$, then $T$ is    ``Lipschitz
for large distances"   in the sense that for each
$\e\ss{0}>0$ there exists  $C=C(\e\ss{0})>0$ so that for all
$\e\ge \e\ss{0}$ and all $x$ in $X$, 
$  T\text{B}_{\e}(x) \subset \text{B}\ss{C\e}(Tx)$.  Similarly, if $T$
is a co-uniformly continuous mapping from a set $X$ onto a convex set
$Y$, then 
$T$ is    ``co-Lipschitz for large distances"  in the sense for each
$\e\ss{0}>0$ there exists  $C=C(\e\ss{0})>0$ so that for all
$\e\ge \e\ss{0}$ and all $x$ in $X$,
$  T\text{B}_{\e}(x)\supset \text{B}\ss{\e/C}(Tx)  $. To see this,
suppose that  $\e$ and $\d$ satisfy $T\text{B}\ss{\e}(x)\supset 
\text{B}\ss{\d}(Tx)$ for every $x$ in $X$.  We have for any $x$ in $X$ that
$$\align
 T\text{B}^X_{{2\e}}(x)&\supset T\left[\cup_{y\in 
\text{B}^X_\e(x)}\text{B}^X_\e(y)\right] = \cup_{y\in \text{B}^X_\e(x)}T
\text{B}^X_\e(y)\\
 &\supset
\cup_{y\in \text{B}^X_\e(x)}\text{B}^Y_\d(Ty)\supset
\cup_{z\in \text{B}^Y_\d(Tx)}\text{B}^Y_\d(z)=\text{B}^Y_{2\d}(Tx),
\endalign$$
 where the last equality follows from the convexity of $Y$.

\proclaim{Proposition 3.4} If the Banach space $Y$ is a uniform
quotient of the Banach space $X$ and $\U$ is a free ultrafilter on
the natural numbers, then $\UY$ is a Lipschitz quotient of $\UX$.
\endproclaim

\demo{Proof} Let $T$ be a uniform quotient mapping from $X$ onto $Y$.  By
the remark,  there is a constant
$C>0$ so that for all
$x$ in $X$ and $r\ge 1$, 
$\text{B}^Y_{Cr}(Tx)\supset T\text{B}^X_{r}(x)\supset 
\text{B}^Y_{r/C}(Tx)$.  This is the only property of $T$ needed in the
proof.  

For each $n$, define $T_n:X\to Y$ by 
$T_nx={{T(nx)}\over n}$.  Then for each $r\ge {1\over n}$ and $x$ in
$X$, we have that
$\text{B}^Y_{Cr}(T_nx)\supset T_n\text{B}^X_{r}(x)\supset 
\text{B}^Y_{r/C}(T_nx)$.  Let $T\ss{\U}:\UX \to \UY$ be the ultraproduct of
the mappings $T_n$, defined for $\tilde{x} =(x_n)$ in $\UX$ by
$T\ss{\U}\tilde{x}=(T_nx_n)$. From the preceding comments it follows
easily that for each $\tilde{x}$ in $\UX$ and $r>0$,
$\text{B}^{\UY}_{Cr}(T\ss{\U}\tilde{x})\supset T\ss{\U}
\text{B}^{\UX}_{r}(x)\supset \text{B}^{\UY}_{r/C}(T\ss{\U}\tilde{x})$, so
that $T\ss{\U}$ is a Lipschitz quotient mapping from $\UX$ onto
$\UY$.
\hfill\qed
\enddemo

We do not know whether a Lipschitz quotient of a separable Banach
space must be a linear quotient of that space.  In the nonseparable
setting there is such an example:  In [AL] it was shown that
$c_0(\aleph)$ is Lipschitz equivalent to a certain subspace $X$ of
$\ell_\infty$, where $\aleph$ is the cardinality of the continuum.  
However, no nonseparable subspace of
$\ell_\infty$  is isomorphic to a quotient of $c_0(\aleph)$ because
this  space (and hence all of its linear quotients) are weakly
compactly generated, while every weakly compact subset of
$\ell_\infty$ is separable.

We now come to the main result we have regarding the linearization of
nonlinear quotients.  Recall that a Banach space $X$ is said to be
{\sl finitely crudely representable\/} in a Banach space $Y$ provided
that there exists $\lambda$ so that every finite dimensional subspace
of $X$ is
$\lambda$-isomorphic to a subspace of $Y$.  If this is true for every
$\lambda>1$, $X$ is said to be {\sl finitely  representable\/} in 
$Y$.

\proclaim{Theorem 3.5} Assume that $X$ is super-reflexive and $Y$  is
a uniform quotient of $X$.  Then $Y^*$ is finitely crudely
representable in $X^*$.  Consequently,  $Y$ is  isomorphic to a
linear quotient of some ultrapower of
$X$. 
\endproclaim

\demo{Proof} By making an arbitrarily small distortion of the norm in $X$,
we can in view of the James-Enflo renorming theorem [DGZ,~p.~149]
assume that $X$ is uniformly smooth.  In view of Proposition~3.4, we
also can assume, via replacement of
$X$ and $Y$ by appropriate ultrapowers, that $Y$ is a Lipschitz
quotient of
$X$.   So for some $\d>0$, there is a mapping $T$ from $X$ onto $Y$
with
$\Lip(T)=1$ and
$T\text{B}_{r}(x)\supset \text{B}_{\d r}(Tx)$ for every $x$ in $X$ and
$r>0$.  Let $E$ be any finite dimensional subspace of $Y^*$,
$E_\perp$ the preannihilator of $E$ in
$Y$,  and 
$Q$ the quotient mapping from $Y$ onto
$Y/E_\perp\equiv E^*$. The composition $QT$ is then a Lipschitz
quotient mapping with $\Lip(QT)\le 1$ and
$QT\text{B}_r(x)\supset \text{B}_{\d r}(QTx)$ for every $x $ in $X$ and
$r>0$.  (Formally, we should here replace $\d$ by an arbitrary
positive number smaller than $\d$, but at the end we know that $Y$ is
reflexive, so the containment we wrote really is true.)  For any
finite dimensional space $Z$, the pair $(X,Z)$ has the UAAP by
Theorem~2.3,  so given any
$\d>\e>0$ there is a ball $B=\text{B}_r(x_0)$ in $X$ and an affine
mapping $G$ from $X$ into $Y/E_\perp$ so that $\sup_{x\in
B}||QTx-Gx||\le \e r$.  This and the quotient property of
$QT$ yield that $G\text{B}_r(x_0)\supset 
\text{B}_{(\d-\e)r}(QT(x_0))$.  Letting
$G_1$ be the linear mapping
$G-G(0)$, we infer that $G\text{B}_1(0)$ contains some ball in
$Y/E_\perp$ of radius $\d-\e$, hence contains the ball around $0$ of
this radius. Therefore
$G_1^*$ is a an isomorphic embedding of $E$ into $X^*$ with
isomorphism constant at most $(\d-\e)^{-1}||G||$.  Since $\Lip(QT)\le
1$, the inequality
$\sup_{x\in B}||QTx-Gx||\le \e r$ implies that $G$ maps the unit ball
of $X$ into some ball of radius at most $1+\e$, so that   $||G||\le
1+\e$.  

From the above we conclude that every finite dimensional subspace of
$Y$ is, for arbitrary $\e>0$,
$(\d^{-1}+\e)$--isomorphic to a subspace of $X^*$, so that $Y^*$ is
finitely crudely representable in  $X^*$.  The ``consequently"
statement is of course a well-known formal consequence of this (since
we now know that $Y$ is reflexive).\hfill\qed
\enddemo

Say that $Y$ is an {\sl isometric Lipschitz quotient\/} of
$X$ provided that for each $\e>0$ there is a mapping $T$ from $X$
onto $Y$ so that $\Lip(T)=1$ and for each $x$ in
$X$, $r>0$, and $0<\d<1$, $T\text{B}_r(x)\supset \text{B}_{\d r}(Tx)$. 
The reason we tracked constants in the proof of Theorem~3.5 was to
make it clear that the following isometric statement is true:

\proclaim{Corollary 3.6}  If $X$ is super-reflexive and $Y$ is an
isometric Lipschitz quotient of $X$, then $Y^*$ is finitely
representable in $X^*$ and
$Y$ is isometrically isomorphic to a linear quotient of some
ultrapower of
$X$.
\endproclaim

We do not know whether a uniform or even Lipschitz quotient of
$\ell_p$,
$1<p\not= 2 < \infty$, must be isomorphic to a linear quotient of
$\ell_p$.  However, a separable  space  which is finitely crudely
representable in
$L_p\equiv L_p[0,1]$ must isomorphically embed into
$L_p$, [LPe].  So we get the following corollary to Theorem~3.5:

\proclaim{Corollary 3.7}  If $Y$ is a uniform quotient of
$L_p$, $1<p<\infty$, then $Y$ is isomorphic to a linear quotient of
$L_p$.  
\endproclaim

Here is another consequence of Theorem~3.5 which is worth mentioning:

\proclaim{Corollary 3.8} A uniform quotient of a Hilbert space is
isomorphic to a Hilbert space.
\endproclaim

Theorem 3.5 and the arguments in [JLS] can be used to classify the
uniform quotients of some spaces other than $L_p$.  As in [JLS],
denote by $\Tt$ a certain one of the 
$2$-convex  modified  Tsirelson-type   spaces, defined as the closed
span of a certain subsequence of the unit vector basis for the
$2$-convexification of the modified Tsirelson space first defined in
[Joh]. (It is known [CO], [CS] that this space is, up to an equivalent
renorming, the space usually denoted by $\Tt$, so our abuse of
notation does little harm.)  From Theorem~3.5 and the results of
[JLS] we deduce the following result, which identifies all the
uniform quotient spaces of a space in a situation where the uniform
quotients differ from the linear quotients.

\proclaim{Corollary 3.9} A Banach space $Y$ is a uniform quotient of
$\Tt$ if and only if $Y$ is isomorphic to a linear quotient of
$\Tt\oplus \ell_2$.
\endproclaim

\demo{Proof} In [JLS] it was proved that $\Tt\oplus \ell_2$ is uniformly
homeomorphic to
$\Tt$,  so every (isomorph  of a) linear quotient of
$\Tt\oplus \ell_2$ is a uniform quotient of $\Tt$. 

Conversely, if $Y$ is a uniform quotient of $\Tt$, then by
Theorem~3.5 $Y$ is isomorphic to a linear quotient of some ultrapower
of $\Tt$.  But by [CS,~p.~150], an ultrapower of $\Tt$ has the form
$\Tt\oplus H$ for some Hilbert space $H$.  This easily implies that
the separable space $Y$ is isomorphic to a linear quotient of
$\Tt\oplus \ell_2$.
\hfill\qed
\enddemo

The next result of this section is the observation that the
$\e$-Fr\'echet derivative of a Lipschitz quotient mapping is a
surjective linear mapping (at least when $\e$ is smaller than the
co-Lipschitz constant of the mapping) and hence the target space is
isomorphic to a linear quotient of the domain space.  
Proposition~3.11 shows that the corresponding statement to
Proposition~3.10 for the G\^ateaux derivative of a Lipschitz quotient
mapping is false.

\proclaim{Proposition 3.10} Suppose $X$, $Y$ are Banach spaces and
$\d>\e>0$. Let
$f : X
\to Y$ be a map which has an $\e$-Fr\'echet derivative, $T$, at some
point
$x_0$ and which satisfies for some $r\ss{0}>0$ and  all
$0<r<r\ss{0}$ the condition $f[\text{\rm B}_r(x_0)]\supset \text{\rm B}_{\d
r}(f(x_0))$.  Then $T$ is surjective.
\endproclaim

\demo{Proof}  Without loss of generality, we may assume that
$x_0=0$, $f(0)=0$, and $\d =1$.  With this normalization we claim
that 
$T\text{\rm B}_1(0)$ contains the interior of $\text{\rm B}_{1-\e}^Y(0)$.  If
not, since $T$ is a bounded linear operator, $\overline{T\text{\rm
B}_{1}(0)}$ cannot contain
$\text{\rm B}_{1-\e}^Y(0)$, and hence there is a vector
$y\ss{0}$ in $\text{\rm B}^Y_1(0)$ so that $d(y\ss{0}, T\text{\rm B}_1(0)) >
\e >0$. (Here, 
$d(\cdot,\cdot)$ denotes the usual infimum ``distance'' between
subsets of
$Y$). Then by linearity, for all $r>0$,
$$ d( ry\ss{0}, T\text{\rm B}_r(0)) > r\e.   \tag3.1 $$

Since $f\text{\rm B}_r(0)\supset \text{\rm B}_{ r}(f(0))$ for all
$0<r<r\ss{0}$, we can choose for $0<r<r\ss{0}$ points $x_{r}$ in
$\text{\rm B}^X_r(0)$ so that 
$f(x_{r})  = ry\ss{0}$. By the definition of   $\e$-Fr\'echet
derivative, we can then choose $0<r\ss{1} < r\ss{0}$ small enough so
that
$$ || f(x_{r}) - Tx_{r} || < r \e  $$  for all $0<r< r\ss{1}$. In
particular, this implies that
$$ d( ry_{0}, T\text{\rm B}_r(0))  \leq ||f( x_{r})  - T x_{r} || < r\e  
 $$ whenever $0<r<r\ss{1}$, contradicting $(3.1)$ above. \hfill\qed
 \enddemo

\proclaim{Proposition 3.11}  There exists a Lipschitz quotient
mapping $f$ from
$\ell_p$, $1\le p<\infty$,  onto itself whose G\^ateaux derivative at
zero is identically zero.
\endproclaim

Define 
$$ f : \left(\sum_0^\infty \  \ell_p\right)_p \to \ell_p
$$ by 
$$f\left(\sum_{n=1}^\infty \sum_{k=0}^\infty a_{nk} e_{nk}\right)= 
\sum_{n=1}^\infty \left(\left[\sum_{k=0}^\infty
g_k(a_{nk}^+)^p\right]^{1\over p} -\left[\sum_{k=0}^\infty
g_k(a_{nk}^-)^p\right]^{1\over p}
\right)e_{n}
$$  where $\{e_{nk}\}_{k=o}^\infty$ is the unit vector basis of the
$n^{\{text{th}}$  copy of $\ell_p$ and 

$$ g_k(t) = \cases
t, &\text{$|t|\ge 2^{-k}$}\\
 0, &\text{$|t|\le
2^{-k-1}$}\endcases$$
and $g_k$ is linear on the intervals
$[2^{-k-1}, 2^{-k}]$ and $[-2^{-k }, -2^{-k-1}]$.

Verifying the Lipschitz condition is easy (check  on the positive cone
and then use  general principles).  

The quotient property is a bit more delicate but not difficult. Here
is the idea:

Suppose that 
$$
\nm{f(x) - y}=d>0;  \quad\quad x= \sum_{n=1}^\infty\sum_{k=0}^\infty
a_{nk}e_{nk}.
$$ We want to find $z$ with $f(z)=y$ and 
$\nm{x-z}\le Cd$ for some constant $C$.  Write
$$ y= f(x) + \sum_{n=1}^\infty b_n e_n
$$ so that 
$$ d^p=\sum _{n=1}^\infty |b_n|^p.
$$

Fix $n$.  We want to  perturb $x$ in the coordinates 
$\enk$ only by a vector $w_n$ with
$\nm{w_n}\le Cb_n$ so that the $n^{\text{th}}$ component of 
$f(x+w_n)$ is the $n^{\text{th}}$ component of $y$.  If we can do this
for each
$n$, then clearly the vector
$z=x+\sum_{n=1}^\infty w_n$  satisfies the requirement.

If $b_n=0$ let $w_n=0$. Otherwise suppose, for definiteness, that
$b_n>0$.  Let 
$A^p = A_n^p = \sum_{k=0}^\infty g_k(a_{nk}^+)^p$. If $A\le b_n$ then
we can let $w_n$ be a multiple, $\alpha$, of $e_{nk_n}$ where $k_n$
is chosen large enough so that
$b_n+a_{nk_n}>2^{-k_n}$.  Then the $n^{\text{th}}$ component of 
$f(x+w_n)$ is
$$
\left(\left[\sum_{k\not= k_n} g_k(a_{nk}^+)^p +
|\alpha+a_{nk_n}|\right]^{1\over p} -\left[\sum_{k\not= k_n}
g_k(a_{nk}^-)^p\right]^{1\over p}
\right)e_{n}.
$$ For $\alpha=10b_n$ this quantity is larger than or equal the 
$n^{\text{th}}$ component of $y$ and for $\alpha\ge0$ small enough
(explicitly,
$\alpha=0$ if $a_{nk_n}>0$ and  $\alpha=|a_{nk_n}|$ otherwise) it is
smaller than the 
$n^{\text{th}}$ component of $y$. By continuity, one can find the
appropriate
$\alpha$.

If $A\ge b_n$, then we let $S=S_n$  be all those $k$'s for which
$a_{nk}> 2^{-k-1}$; that is, for which 
$g_k(a_{nk})>0$.  In this case we can let $w_n$ be a multiple of
$\sum_{k\in S} a_{nk}^+ e_{nk}$.
\hfill\qed

In order to ``soup-up" the example in Proposition~3.11, we need the
following perturbation lemma for co-Lipschitz mappings:

\proclaim{Lemma 3.12}  Suppose that $f$ and $g$ are continuous
mappings from the Banach space $X$ into the Banach space  $Y$ and 
$\text{Lip}(g)<\text{co-Lip}(f)^{-1}$.  Then $f+g$ is co-Lipschitz
and  
$\text{co-Lip}(f+g)\le [{1-\text{co-Lip}(f) \text{Lip}(g)}]^{-1}$.
\endproclaim

\demo{Proof} Lemma 3.12 follows from a classical   successive approximations
argument.  By dividing $f$ and $g$ by $\text{co-Lip}(f)$, we can
assume that $1=\text{co-Lip}(f)> \text{Lip}(g) \equiv \d $.  Let $x$
be in $X$ and $r>0$.  Given $y\in Y$ with $||y||< r$, we need to find
$z\in X$ with $||z||< (1-\d)^{-1}r$ and $f(x+z)+g(x+z)=f(x)+g(x)+y$.  
Set
$z_0=0$ in $X$.   Recursively choose $z_k$ in $X$  so that for each
$n=0,1,2,\dots$,  $||z_n|| <\d^{n-1} r$   and 
$$ f(x+\sum_{k=0}^{n+1} z_k) + g(x+\sum_{k=0}^n z_k)= f(x)+g(x)+y. 
\tag3.2
$$ We can make the choice of $z_1$ because $\text{co-Lip}(f)=1$.  If
$(3.2)$ holds for $n$ with $||z_{n+1}||<\d^n r$, then we have 
$||g(x+\sum_{k=0}^{n+1} z_k)-g(x+\sum_{k=0}^n z_k)||<\d^{n+1}r$. 
Again using the condition $\text{co-Lip}(f)=1$, we can choose
$z_{n+2}$ with 
$||z_{n+2}||<\d^{n+1} r$ so that $(3.2)$ holds with $n$ replaced by
$n+1$.  Setting 
$z=\sum_{k=1}^\infty z_k$, we get  the desired conclusion because $f$
is continuous.\hfill\qed
\enddemo

\proclaim{Corollary 3.13} Let $T$ be any bounded linear operator on
$\ell_p$, $1\le p <\infty$.  Then there is a Lipschitz quotient
mapping from $\ell_p$ onto itself whose G\^ateaux derivative at zero
is $T$.
\endproclaim

\demo{Proof} Add $T$ to a suitable multiple of the $f$ from
Proposition~3.11. \hfill\qed
\enddemo

If the point of G\^ateaux differentiability is an isolated point in
its level set, the phenomenon in Proposition~3.11 cannot occur:

\proclaim{Proposition 3.14} Suppose that $f$ is a Lipschitz quotient
mapping from
$X$ to $Y$, $f$ has G\^ateaux derivative $T$ at some point
$p$ in $X$, and $p$ is isolated in the level set $[f=f(p)]$.  Then
$T$ is an isomorphism from $X$ into $Y$.
\endproclaim

\demo{Proof}  Without loss of generality, we may assume that $p=0,$
$f(0)=0$, and
$f\text{\rm B}_r(x)\supset \text{\rm B}_r(fx)$ for each $r>0$ and $x$ in $X$.
For 
$t>0$, define
$f_{t}: X \to Y$ by 
$$ f_{t}(x) = {{f(tx)}\over{t}}. $$ The $f_{t}$ then converge
pointwise as
$t\to 0$  to the G\^ateaux derivative
$T$,  and moreover $f_t\text{\rm B}_r(x)\supset \text{\rm B}_r(f_tx)$ for each
$r>0$, $t>0$,  and $x$ in $X$.

If $T$ is not an isomorphism, then there exists a unit vector
$x \in X$ such that $||Tx||<{1\over 4}$ and hence $||f_{t}x|| <1/4$
for all
$0<t< t\ss{0}$. This implies that for each $0<t< t\ss{0}$, there
exists
$x_{t}$ in $X$ such that 
$|| x_{t} - x  || \leq 1/4$ and $f_{t} x_{t}  =0$. Hence,
$f(tx_{t}) =0$, and  the family of points
$\{ tx_{t} \} \subset [f=0]$ tends to zero as $t \to 0$.
\hfill\qed
\enddemo
 
If one composes the example in Proposition~3.11 with the projection
onto the first coordinate of $\ell_p$, then one obtains a Lipschitz
quotient mapping-- call it $g$--from $\ell_p$ onto the real line
which has zero G\^ateaux derivative at
$0$.  This phenomenon of course cannot happen when the domain space
is finite dimensional, since then the G\^ateaux derivative is a
Fr\'echet derivative. In particular, the restriction of $g$ to any
finite dimensional subspace of
$\ell_p$ is not a Lipschitz quotient mapping.  However, a nonlinear
quotient mapping onto a separable space does have a separable
``pullback":

\proclaim{Proposition 3.15} Let $f$ be a continuous, co-Lipschitz
(respectively, co-uniformly continuous)  mapping from the  metric
space $X$ onto the  separable metric space
$Y$ and let
$X\ss{0}$ be a separable subset of
$X$.  Then there is a separable closed subset $X\ss{1}$ of $X$ which
contains $X\ss{0}$ so that the restriction of $f$ to $X\ss{1}$  is 
co-Lipschitz (respectively, co-uniformly continuous), with
co-Lip$(f_{|X_1})\le$co-Lip$(f)$ when $f$ is co-Lipschitz.  If $X$ is
a Banach space, $X\ss{1}$ can be taken to be a subspace of $X$.
\endproclaim

\demo{Proof} For definiteness, assume that $f$ is co-Lipschitz, normalized
so that $\text{co-Lip}(f)=1$, so that for each
$x$ in $X$ and $r>0$,
$f[\text{\rm B}^X_r(x)]\supset \text{ int}\left(\text{\rm
B}^Y_r(f(x))\right)$.  Let
$W\ss{0}$ be a countable dense subset of $X\ss{0}$ and build countable
subsets 
$W\ss{0} \subset W\ss{1}\subset W\ss{2} \subset  \dots$ of $X$ so that
for each $n$,  each $x$ in
$W\ss{n}$, and each positive rational number $r$,
$$
\overline{f\left[\text{\rm B}^X_r(x)\cap W\ss{n+1}\right]} \supset \text{\rm
B}^Y_r(f(x)).
\tag3.3
$$
 When $X$ is a Banach space, we can also make sure that $W_n$ is
closed under rational linear combinations for $n\ge 1$. 

Note that  $(3.3)$ must hold for  all positive
$r$.   Let $X\ss{1}$ be the closure of $\cup_{n=0}^\infty W\ss{n}$. 
From
$(3.3)$ we deduce that for each $x$ in $X\ss{1}$ and each $r>0$,
$$
\overline{f\left[\text{\rm B}^{X_1}_r(x)\right]} \supset \text{\rm B}^Y_r(f(x)).
\tag3.4
$$ 
We complete the proof by observing that in $(3.4)$ we can remove
the closure if we replace the right hand side by its interior.  For
linear $f$ this is sometimes called the ``little open mapping
theorem"; in fact, the
 successive approximation argument requires in addition to $(3.4)$
only continuity of
$f$. (If
$||y-f(x\ss{0})|| < r-\tau$, choose $x_1$ so that
$||x_1-x_0||<r-\tau$ and
$||y-f(x_1)||<{\tau \over 2}$.  Then choose $x_n$ recursively to
satisfy $||x_n-x_{n-1}||<{\tau\over{2^{n-1}}}$ and
$||y-f(x_n)||<{\tau\over{2^n}}$.  Clearly $\xn$ converges to some
point, $x$, in $\text{\rm B}^{X_1}_r(x_0)$, and $f(x)=y$ by the continuity
of
$f$.) \hfill\qed
\enddemo

{\bf Remark 3.16.} Note that the last part of the argument
for Proposition~3.15 shows that if $f$ is a continuous mapping from
the metric space $X$ onto a metric space $Y$, $X_0\subset X$, and the
restriction of $f$ to $X_0$ is co-Lipschitz (respectively,
co-uniformly continuous) when considered as a mapping onto $f[X_0]$,
then $f$ maps the closure of $X_0$ onto the closure of $f[X_0]$ and
$f_{|\overline{X_0}}$ is co-Lipschitz (respectively, co-uniformly
continuous) when considered as a mapping onto $\overline{f[X_0]}$.

\proclaim{Corollary 3.17} Suppose that $f$ is a continuous,
co-Lipschitz (respectively, co-uniformly continuous)  mapping from the
Banach space $X$ onto the Banach space $Y$, $X_0$ is a separable
subspace of $X$, and $Y_0$ is a separable subspace of $Y$.  Then there
exist separable closed subspaces $X_1$ of $X$, $Y_1$ of $Y$ so that
$X_0\subset X_1$, $Y_0\subset Y_1$, $f[X_1]=Y_1$, and the restriction
of
$f$ to $X_1$ is a co-Lipschitz (respectively, co-uniformly continuous)
mapping of $X_1$   onto $Y_1$.
\endproclaim

\demo{Proof} We prove the co-Lipschitz case; the co-uniformly continuous
case is similar. Note that if $Z$ is a subset of $Y$, then the
restriction of
$f$ to
$f^{-1}[Z]$ is co-Lipschitz when considered as a mapping onto $Z$,
with co-Lip$(f_{|f^{-1}[Z]})\le $ co-Lip$(f)$. Using this and
Proposition~3.15, we build separable closed subsets $X_0\subset
W_0\subset W_1\subset W_2\subset \cdots$ so that for each $n$, the
restriction of $f$ to $W_n$ is a co-Lipschitz mapping onto $f[W_n]$
with co-Lip$(f_{|W_n})\le \text{co-Lip}(f)$, 
$f[W_0]\supset Y_0$, $W_{n+1}\supset \text{\rm span}\, W_n$, and
$f[W_{n+1}]\supset {\text{\rm span}}\,Y_n$. From this it follows that the
restriction of
$f$ to the separable (possibly nonclosed) subspace $\cup_{n=0}^\infty
W_n$ is a co-Lipschitz mapping onto its image, which is a (possibly
nonclosed) subspace.  The desired conclusion now follows from
Remark~3.16.\hfill\qed
\enddemo

Our final result of this section, while modest, seems to be new even
for bi-Lipschitz equivalences.  Note that it yields that $c_0$ is not 
bi-Lipschitz equivalent to $C[0,1]$, a result first proved in [JLS].

\proclaim{Theorem~3.18} If $X$ is  Asplund  and $Y$ is a Lipschitz
quotient of $X$, then $Y$ is Asplund.
\endproclaim

\demo{Proof} Recall (see e.g. Theorem~5.7 in [DGZ]) that a Banach space is
Asplund if and only if every separable subspace has separable dual. 
Therefore, in view of Proposition~3.17, we may assume that $X$ is
separable, and need to prove that $Y^*$ is separable.  Let $f$ be a
Lipschitz quotient mapping from $X$ onto $Y$ with $\Lip(f)=1$, set
$C:=\text{ co-Lip}(f)$, and assume that
$Y^*$ is nonseparable.  Then there is an uncountable set
$\{y_\gam^*\}_{\gam\in
\Gam}$  of norm one functionals in $Y^*$ so that for each
$\gam\not=\gam'$, $||y_\gam^*-y_{\gam'}^*|| > 1/2$. Set
$f_\gam=y_\gam^* f$;  then Lip$(f_\gam)\le \text{\rm Lip}(f)=1$.   By
[Pre], each
$f_\gam$ has a Fr\'echet derivative at some point $x_\gam$ in the
unit ball of
$X$. 
\enddemo

Let $\e=(10C)^{-1}$ and for each $\gam$ choose $\d_\gam>0$ so that
whenever
$||z||\le \d_\gam$,
$$ |f_\gam(x_\gam+z)-f_\gam(x_\gam)-f_\gam'(x_\gam;z)| \le
\e||z||.\tag3.5
$$
 By passing to a suitable uncountable subset of $\Gam$, we may
assume that
$\d:=\inf_{\gam\in\Gam} \d_\gam >0$, and also (since $X$ and $X^*$ 
are separable)  that for all
$\gam$,
$\gam'$
$$||x_\gam-x_\gam'||<\e \d, \tag3.6
$$
$$ ||f_\gam'(x_\gam)-f_{\gam'}'(x_{\gam'})||<\e, \tag3.7
$$
$$ |f_\gam(x_\gam)-f_{\gam'}(x_{\gam'})|<\e\d.\tag3.8
$$ 
If $||z||\le \d$, we have:
$$\align
|f_\gam(x_\gam & +z)- f_{\gam'}(x_\gam+z)| = 
|\left(f_\gam(x_\gam+z)-f_\gam(x_\gam)-f_\gam'(x_\gam;z)\right)+
\left(f_{\gam}(x_{\gam})-f_{\gam'}(x_{\gam'})\right)
\\
& +
\left(f_\gam'(x_\gam;z)-f_{\gam'}'(x_{\gam'};z)\right)+
\left(-f_{\gam'}(x_{\gam'}+z)+f_{\gam'}(x_{\gam'})+f_{\gam'}'(x_{\gam'};z)\right
)\\
 &+
\left(f_{\gam'}(x_{\gam'}+z)-f_{\gam'}(x_\gam+z)\right)|\\
 &\le \e||z||+\e\d+||f_\gam'(x_\gam)-f_{\gam'}'(x_{\gam'})|| 
||z||+
\e||z|+||x_\gam-x_{\gam'}|| \ \,  \text{by $(3.5)$, $(3.8)$, $(3.5)$}\\
 & <
5\e\d=\d/(2C).\qquad\qquad\qquad\qquad\qquad\qquad\qquad\qquad\qquad\
\ \ \ 
\text{by  $(3.7)$, $(3.6)$}
\endalign$$ 
On the other hand, 
$$\align
\sup_{z\in \text{\rm B}_\d^X(0)} |f_\gam(x_\gam  +z)- f_{\gam'}(x_\gam+z)|
& =
\sup_{z\in \text{\rm B}_\d^X(0)}
|\left(y_\gam^*-y_{\gam'}^*\right)f(x_\gam+z)|  \\
 & \ge
\sup_{y\in \text{\rm B}_{\d/C}^Y(0)} 
|\left(y_\gam^*-y_{\gam'}^*\right)(f(x_\gam)+y)| 
\cr & \ge ||y_\gam^*-y_{\gam'}^*||\left(\d/C\right) >
\d/(2C).\qquad\qquad\qquad\qquad\ \ \text{\qed}
\endalign$$

{\bf Remark 3.19.}  The analogue of Theorem~3.18 for uniform
quotient mappings is false.  Ribe [Rib2] gave an example of a
separable, reflexive space which is bi-uniformly homeomorphic to a
space which contains an isomorphic copy of $\ell_1$.

\head 4. Examples related to the Gorelik Principle\endhead

One natural way of constructing a uniform or Lipschitz quotient
mapping is to follow a bi-uniformly continuous or bi-Lipschitz
mapping with a linear quotient mapping.  By the Gorelik principle
[Gor], [JLS], the resulting mapping cannot map a ``large" ball in a
finite codimensional subspace of the domain into a ``small" ball
(Theorem~1.1 in [JLS] gives a precise quantitative meaning to this
statement).  It is natural to guess that any uniform quotient mapping
satisfies this version of the Gorelik principle.  If this were true,
this would provide the machinery to prove, for example, that every
uniform quotient of
$\ell_p$, $1<p<\infty$, is isomorphic to a linear quotient of
$\ell_p$.  Unfortunately, Proposition~4.1 shows that there is not a
Gorelik principle for uniform quotient mappings.
 
\proclaim{Proposition 4.1} Let $X$ be a Banach space and set 
$Z=X\oplus_1 X\oplus_1 \R$. Then there exists a Lipschitz, co-uniform
mapping
$T$ from $Z$ onto $X$ so that $T\left(\text{\rm B}_1^{X\oplus
X}(0)\right)=\{0\}$.
\endproclaim

\demo{Proof} We shall define $T(x,y,\lam)= ag(x,\lam)+f(y,\lam)$ for
appropriate Lipschitz functions $f$ and $g$ from
$X\oplus\R$ into $X$ and for an appropriately small $a$; namely, for
$a={1\over 2}$. For each integer $k$ let $\ukj$ be a maximal $4\cdot
2^k$-separated set in $X$, and define for each $y$ in $X$
$$ f(y,2^k)= \sum_{j=1}^\infty \left(2^k-\left| \,
\nm{y-u\ss{k,j}}-2^k\right| \,  \right)^+
{{y-u\ss{k,j}}\over\nm{y-u\ss{k,j}}}
$$ (where ${0\over 0} \equiv 0$.)  Write $f_k(y)=f(y,2^k)$; the
formula for $f$ just says that  $f_k$ translates the ball in $X$ of
radius
$2^k$ around each $u\ss{k,j}$ to a ball around the origin, that $f_k$
vanishes on the complement of $\dstyle
\cup_{j=1}^\infty \text{\rm B}_{2^{k+1}}(u\ss{k,j})$, and that $f_k$ is
affine on the intersection  of $\text{\rm B}_{2^{k+1}}(u\ss{k,j}) \sim
\text{\rm B}_{2^{k}}(u\ss{k,j})$ with any ray emanating from $u\ss{k,j}$.  
Extend
$f$ to a function from $X\oplus \R$ into $X$ by making $f$ affine on
each interval of the form $[(y,2^k), (y,2^{k+1})]$, by setting
$f(y,0)=0$, and by defining $f(y,\lam)=f(y,-\lam)$ when
$\lam<0$.  It is apparent that each $f_k$ has Lipschitz constant
one,  that $f$ has Lipschitz constant at most two, and that
$||f(y,\lam)|| \le |\lam|$.
\enddemo

Define $g$ by $g(x,\lam)=\left(\left[ |\lam| \vee \left(
||x||-1\right)\right]\wedge 1 \right)  x$. It is easy to see that $g$
has Lipschitz constant one and hence $T$ has Lipschitz constant at
most three (as long as $a\le 1$).  Evidently $T$ vanishes on the unit
ball of the hyperplane
$[\lam=0]$, so it remains to check that $T$ satisfies the uniform
quotient property.

Let $(x\ss{0}, y\ss{0}, \lam\ss{0})$ be in $Z$ and $r>0$.  We need to
check that $T\text{\rm B}_{r}(x\ss{0}, y\ss{0}, \lam\ss{0})$ contains a
ball around
$T(x\ss{0}, y\ss{0}, \lam\ss{0})$ whose radius can be estimated from
below in terms of $r$ only.  We can assume, without loss of
generality, that $r<1$.

{\bf Case 1:}  $|\lam\ss{0}|<{r \over {48}}$.  

This is the only case in which the mapping $f$ comes into play and
also the only case in which we need to consider points off the
hyperplane $[\lam=\lam\ss{0}]$ in order to verify the quotient
property.  We show in this case that $T\text{\rm B}_{r}(x\ss{0}, y\ss{0},
\lam\ss{0})$ contains a ball around
$T(x\ss{0}, y\ss{0}, \lam\ss{0})$ of radius $r\over {32}$.  By
symmetry we can assume that $\lam\ss{0}\ge 0$.

Define $k$ by $2^k \le {{r}\over 6}<2^{k+1}$ (so 
$0<2^k-\lam\ss{0}\le  2^k$) and choose $j$ to satisfy  
$||y\ss{0}-u\ss{k,j}||<2^{k+2}$.  Then 
$f_k\text{\rm B}_{2^k}(u\ss{k,j}) = \text{\rm B}_{2^k}(0) $  and 
$\{x\ss{0}\}\oplus\text{\rm B}^X_{2^k}(u\ss{k,j})\oplus \{2^k\}\subset 
\text{\rm B}_{6\cdot 2^k}(x\ss{0},y\ss{0},\lam\ss{0})\subset \text{\rm
B}_r(x\ss{0}, y\ss{0}, \lam\ss{0})$.  Therefore 
$$\align \text{\rm B}_r(x\ss{0}, y\ss{0}, \lam\ss{0}) & \supset 
ag(x\ss{0}, 2^k) + f_k \text{\rm B}_{2^k}(u\ss{k,j}) \supset ag(x\ss{0},
\lam\ss{0})+\text{\rm B}^X_{2^k-a(2^k-\lam\ss{0})}(0) \\
 &\supset ag(x\ss{0}, \lam\ss{0})+f(y\ss{0},\lam\ss{0})+ \text{\rm
B}^X_{2^{k-1}-(1-a)\lam\ss{0}}(0)
\supset \text{\rm B}^X_{r\over 32}\left(T(x\ss{0}, y\ss{0},
\lam\ss{0})\right).
\endalign$$

{\bf Case 2:} $|\lam\ss{0}|\ge {r\over 48}$.

Again we assume by symmetry that $\lam\ss{0}\ge 0$.  Fix $(x\ss{0},
y\ss{0}, \lam\ss{0})$ in $Z$ and set 
$s\ss{0}=\left[ \lam\ss{0} \vee \left(
||x\ss{0}||-1\right)\right]\wedge 1$  so that
$g(x\ss{0},\lam\ss{0})=s\ss{0} x\ss{0}$ and $r\le 48(\lam\ss{0}\wedge
1)\le 48 s\ss{0}$.  Fix $z$ in $X$ with
$||z||\le {r\over 5}$.  We want to find a vector $x$ in $X$ with
$||x-x\ss{0}||\le r$ and
$g(x,\lam\ss{0})=s\ss{0}(x\ss{0}+z) \  (=g(x\ss{0},\lam\ss{0}) +
s\ss{0} z)$.  This will give
$$\align T\text{\rm B}_r(x\ss{0}, y\ss{0}, \lam\ss{0}) & \supset 
f(y\ss{0},  \lam\ss{0}) + a g\left[T\text{\rm B}^X_r(x\ss{0})\oplus
\{\lam\ss{0}\}\right] \\
 & \supset  f(y\ss{0},  \lam\ss{0}) + a
g(x\ss{0},\lam\ss{0}) +  a s\ss{0} \text{\rm B}^X_{r\over 5}(0) \supset
\text{\rm B}^X_{{a r s_0 }\over 5}\left(T(x\ss{0}, y\ss{0},
\lam\ss{0})\right).
\endalign$$
 Since ${{a r s_{0} }\over 5}\ge {{ar^2}\over 240}$, this will
verify that $T$ has the uniform quotient property.

The vector $x$ will be of the form $x=t(x\ss{0} + z)$ for appropriate
$t$ close to one.  Write 
$s(t) = \left[\lam\ss{0} \vee \left(t||x\ss{0}+z||-1\right)\right]
\wedge 1$, so that 
$g(t(x\ss{0}+z),\lam\ss{0})=s(t)t(x\ss{0}+z)$.  We need to find
$t\ss{0}$ so that $t\ss{0}s(t\ss{0})=s\ss{0}$ and 
$||x\ss{0}-t\ss{0}(x\ss{0}+z)|| \le r$.  If $||x\ss{0}||\ge 3$, the
choice $t\ss{0}=1$ works, since then $||x\ss{0}+z||-1\ge 1$ and thus
$s(1)=s\ss{0}$.  so assume that $||x\ss{0}||\le 3$.  In this case, by
continuity of $s(\cdot)$, it is enough to check that 
$(1-{r\over 4})s(1-{r\over 4})\le s\ss{0} \le (1+{r\over
4})s(1+{r\over 4})$; actually, we check the stronger inequalities:
$$ s(1-{r\over 4})\le s\ss{0}
\leqno\text{\rm (L)}
$$

$$ s\ss{0} \le s(1+{r\over 4}).
\leqno\text{\rm (R)}
$$ We first check (R).  If $s\ss{0}=\lam\ss{0}\wedge 1$, then $s(t)\ge
s\ss{0}$ for all $t\ge 0$ and so (R) is clear.  Otherwise 
$s\ss{0}=||x\ss{0}||-1$ and we get 
$$\align
 & (1+{r\over 4})||x\ss{0}+z||-1 \ge (1+{r\over
4})(1+s\ss{0})-(1+{r\over 4}){r\over 5} -1 \\
\ge & (1+{r\over 4})s\ss{0}+r\left({1\over 4}-{{4+r}\over {20}}\right)
\ge (1+{r\over 4})s\ss{0}.
\endalign$$ 
To check (L), suppose first that $s\ss{0}\ge ||x\ss{0}||-1$. Then
$$ (1-{r\over 4})||x\ss{0}+z||-1 \le  (1-{r\over 4}) (1+s\ss{0}) +
(1-{r\over 4}) {r\over 5} -1
\le (1-{r\over 4})s\ss{0},
$$ and (L) follows.  On the other hand, if 
$||x\ss{0}||-1> s\ss{0}$, then $||x\ss{0}||>1$, and hence
$$ (1-{r\over 4})||x\ss{0}+z||-1 \le ||x\ss{0}||-1-{r\over 4} +
(1-{r\over 4}) {r\over 5}< ||x\ss{0}||-1,
$$ which also yields (L). \hfill\qed

Of course, if $X$ is isomorphic to its square and also to its
hyperplanes (for example, if $X$ is $\ell_p$ or $L_p$, $1\le p\le
\infty$), then Proposition 4.1 says that there is a uniform quotient
mapping from
$X$ onto itself which maps the unit ball of a hyperplane to zero.

Proposition 4.1 also implies that for each $n$ there is a uniform
quotient mapping from $\Rn{2n+1}$ onto $\Rn{n}$ which maps the unit
ball of $\Rn{2n}$ to zero.  However, in the finite dimensional case,
more can be said:

\proclaim{Proposition 4.2} There is a Lipschitz map $T$ from 
 $\Rn{3}=\Rn{2}\oplus\R$ onto $\Rn{2}$ such that $T$ is a co-uniform
quotient map and
 $T\text{\rm B}_1^{\Rn{2}\oplus\{0\}}(0)=0$.
 \endproclaim

 \demo{Proof}. For $\theta\in\R$ let
$$
 U_\theta=\left(\matrix
 \cos\theta&-\sin\theta\\
                       \sin\theta&\cos\theta\endmatrix\right).
$$
 For $0\le a\le1$ let $r_a:\R^{+}\to[0,1]$ be defined by

$$ r_a(t)=\cases
a,    &\text{if $0\le t \le 1$;}\\
         1-(1-a)(2-t), &\text{if $1<t<2$;}\\
           1,          &\text{if $ 2<t$}.\endcases
$$
 We also let $r_a(t)=1$ if $a>1$ and $r_a(t)=r_{|a|}(t)$, and define
$T:\Rn{3}\to\Rn{2}$ by
$$ T(x,a)=r_a^2(\nm{x})U_{\theta(r_a(\nm{x}))} x,$$ where
$\theta:\R^{+}\to\R^{+}$ is defined by
$$
\theta(t)=2\pi/t.
$$
\enddemo

We check first that $T$ is Lipschitz. It is clearly enough to show
that its restriction to the set $\{(x,a): \nm{x}\le 2\}$  is
Lipschitz,   which follows immediately by noting that  it is the
composition of 
$$(x,a)\in \{(x,a): \nm{x}\le 2\}\to (x,r_a(\nm{x})),$$ which is
clearly Lipschitz, followed by
$$(x,t)\in \{(x,t): \nm{x}\le 2, 0\le t\le 1\}
\to t^2 U_{\theta(t)} x,$$ which has bounded partial derivatives.

It remains to show that $T$ is co-uniform.  We note first that for
each
$a>0$ the function 
$$f_a(x)=T(x,a)$$  is a Lipschitz homeomorphism of the plane, with
inverse given by
$$ f^{-1}_a(y)=(s_a(\nm{y})/\nm{y})U_{-\theta(r_a(s_a(\nm{y})))} y,$$
where
$s_a$ is inverse of $t\to tr^2_a(t)$.  The function $tr^2_a(t)$ has a
positive derivative bounded below by $a^2\wedge 1$; moreover, if
$0<\tau<1$, this derivative is bounded below by $\tau^2$ on the
interval $t>1+\tau$. Thus,  if either $a\ge\tau>0$ or $\nm{x}\ge
1+\tau$,  the mapping
$f^{-1}_a$ has its derivative at $f_a(x)$ bounded in norm by a
constant
$\kappa_\tau$ depending only on $\tau$. Consequently,
$$f_a B_\tau(x)\supset f_a B_{\tau/\kappa_\tau}(f_a(x)),$$ if either
$a\ge\tau$ or $\nm{x}\ge 1+\tau$.

To check co-uniformity of $T$, take any $(x,a)\in\Rn{3}$ with $a>0$
and let $0<r<1$.

{\bf Case 1:} $r\le 10 a$ or $\nm{x}>1+r/10$. In this case the above
inclusions show that
$$T B_r(x,a)\supset f_a  B_r(x)\supset B_{r\kappa_{r/10}}(f_a(x)) =
B_{r\kappa_{r/10}}(T(x,a)).$$

{\bf Case 2:}  $r>10 a$ and $\nm{x}\le 1$.  Note first that by the
definition of $f_a$,
$$\{f_\gamma(\beta^2x/\gamma^2): 1/(k+1)<\gamma\le 1/k\} =\{y:
\nm{y}=\beta^2\nm{x}\},$$ whenever $k$ is an integer and $\beta\le
1/(k+1)$. Let now $k_0$ be the integer so that 
$1/(k_0+1)< r \le 1/k_0$. Let $k\ge 3k_0$, and let 
$1/(k+2)<\beta\le 1/(k+1)$ and $1/(k+1)<\gamma\le 1/k$. We have
$$\align 
\nm{(x,a)-(\beta^2x/\gamma^2,\gamma)}^2  &\le
(1-\beta^2/\gamma^2)^2+(a-\gamma)^2
\le (1-(k/(k+2))^2)^2+1/(3k_0)^2\\
 &\le 4/(3k_0+2)^2+1/(3k_0)^2
\le r^2.
\endalign$$ 
Hence $T B_r(x,a)$ contains all the points 
$f_\gamma(\beta^2x/\gamma^2)$ with $\beta$ and $\gamma$ as above.
Consequently,
$$T B_r(x,a)\supset B_{\nm{x}/(3k_0+1)^2}(0)\supset
B_{\nm{x}r^2/16}(0).$$ Since $\nm{T(x,a)}=a^2\nm{x}\le \nm{x}r^2/100$,
we conclude that
$$T B_r(x,a)\supset B_{\nm{x}r^2/50}(T(x,a)).$$ If $\nm{x}\ge r/2$,
this means that 
$$T B_r(x,a)\supset B_{r^3/200}(T(x,a)),$$ while if $\nm{x}< r/2$,
then
$$T B_r(x,a)\supset T B_{r/2}(rx/(2\nm{x}),a)\supset 
B_{(r/2)(r/2)^2/50}(T(x,a))\supset  B_{r^3/400}(T(x,a)).$$

{\bf Case 3:}  $r>10 a$ and $1<\nm{x}\le 1+r/10$.  Put $u=x/\nm{x}$.
Then $B_r(x,a)\supset B_{9r/10}(u,a)$ and by the proof of case~2,
$$T B_{9r/10}(u,a)\supset B_{(9r/10)^2/16}(0)\supset B_{r^2/20}(0).$$
Also $r_a(\nm{x})=1-(1-a)(2-\nm{x})\le a+\nm{x}-1\le r/5$. Thus
$\nm{f_a(x)}\le 11r^2/250$, and therefore
$$T B_r(x,a)\supset B_{r^2/400}(T(x,a)).$$

\hfill\qed

We do not know whether there is a Gorelik principle for Lipschitz
quotient mappings.  We do not
know, for example, if there is a Lipschitz quotient mapping from
some space of dimension at least three onto
$\Rn{2}$ which vanishes on a hyperplane.  Notice that if
$f:X\to Y$ is a Lipschitz quotient mapping which sends the unit ball
of some subspace $Z$ of $X$ to zero, and $\U$ is a free ultrafilter
on the positive integers, then there is a Lipschitz quotient mapping
$\Uf$ from the ultrapower
$\UX$ of
$X$ onto the ultrapower $\UY$ of $Y$ which sends the entire {\sl
subspace\/} $\UZ$ of $\UX$ to zero.  Indeed, define $f_n:X\to Y$ by
$f_n(x)=nf\left({x\over n}\right)$ and let $\Uf$ be the ultraproduct
of the $f_n$'s, defined by 
$\Uf(x_1,x_2,\dots)=(f_1(x_1),f_2(x_2),\dots)$.  Now if $X$ is finite
dimensional, then so is $Y$, and $X=\UX$, $Y=\UY$, and
$Z=\UZ$, so one obtains a Lipschitz quotient mapping from $X$ onto
$Y$ which maps the subspace $Z$ to zero.  (Incidentally, if $f$ is a
uniform quotient mapping from a space $X$ onto a space $Y$ which maps
a subspace $Z$ of $X$ to zero, then the construction of
Proposition~3.4 produces a {\sl Lipschitz\/} quotient mapping from
$\UX$ onto $\UY$ which maps $\UZ$ to zero.)

In the case of mappings from $\Rn{n}$ to $\Rn{n}$, there is a close
connection between Lipschitz quotient mappings and quasiregular
mappings ([Ric] is the standard source for
 quasiregular mappings).  Recall that a map $f$ from a domain
$G$ in $\Rn{n}$ to $\Rn{n}$ is called quasiregular
provided\hfill\break (i) $f$ is $ACL^n$, i.e., $f$ is continuous, its
restriction to every line in the direction of each of the coordinate
axes is absolutely continuous and its partial derivatives belong
locally to the space $L_n(\Rn{n})$.
\hfill\break (ii) The $n$ by $n$ matrix $D(x)$ of partial derivatives
of $f$ satisfies
$\|D(x)\|^n\le KJ(x)$ for almost every $x$, where $K$ is a constant,
$J(x)$ is the determinant of $D(x)$, and $\|D(x)\|$ is its norm as an
operator from
$\ell_2^n$ to itself.

If $f: \Rn{n}\to \Rn{n}$ is a Lipschitz quotient mapping then (i)
holds  trivially. As for (2), we have the somewhat weaker
statement\hfill\break (ii') For every $x$ at which $f$ is
differentiable
$$
\|D(x)\|^n\le K|J(x)|
$$ where $K=(M/m)^{n-1}$ with $M$ (respectively, $m$)   the Lipschitz 
(respectively, co-Lipschitz) constant of $f$.

It seems likely that  some results from the quasiregular theory can
be carried over to the case of Lipschitz quotient mappings. A deep
result of Reshetnyak (which is nicely presented in [Ric]) states that
the level sets of a nonconstant quasiregular mapping from a domain
$G$  in $\Rn{n}$ into
$\Rn{n}$ are discrete sets.  In Proposition~4.3 we prove that the
level sets of a co-Lipschitz continuous mapping from $\Rn{2}$ to
$\Rn{2}$ are discrete.  Note that  this result implies in particular
that a Lipschitz quotient mapping from
$\Rn{2}$ to
$\Rn{2}$ cannot vanish on an interval, which means that some version
of the Gorelik principle is true for such mappings.

\proclaim{Proposition 4.3} Let $f:\Rn{2}\to\Rn{2}$ be a continuous
and  co-Lipschitz mapping. Then for every $y\in \Rn{2}$ the set
$f^{-1}(y)$ is  discrete.
\endproclaim

\demo{Proof} We use the following  simple lemma concerning the lifting of
Lipschitz curves:
\enddemo

\proclaim{Lemma 4.4} Suppose that $f:\Rn{n}\to X$ is  continuous and
co-Lipschitz with constant one, $f(x)=y$, and $\xi:[0,\infty) \to X$
is a curve with Lipschitz constant one, and $\xi(0)=y$.  Then there
is a curve $\phi:[0,\infty)\to \R^{n}$ with Lipschitz constant one
such that $\phi(0)=x$ and $f(\phi(t))=\xi(t)$ for $t\ge 0$.
\endproclaim

\demo{Proof}  For $m=1,2,\ldots$ define $\phi_{m}(0)=x$, and, by induction, 
assuming that $f(\phi_{m}({{k}\over{m}})) = \xi ({{k}\over {m}})$,
choose $\phi_{m}({{k+1}\over {m}})$ such that
$\|\phi_{m}({{k+1}\over {m}})-\phi_{m}({{k }\over {m}})\|\le
{{1}\over {m}}$
 and $f(\phi_{m}({{k+1}\over {m}}))= \xi({{k+1}\over {m}})$. Extend
$\phi_{m}(t)$ to a Lipschitz curve
$\phi_{m}:[0,\infty)\to \R^{m}$ having Lipschitz constant one. The
limit $\phi$ of any convergent subsequence of
$\phi_{m}$ has the desired properties. \hfill 
\qed
\enddemo

Without loss of generality, assume $\text{\rm B}_r(f(x))\subset  f(\text{\rm
B}_r(x ))$ for every $x$ in $\Rn{2}$ and
$r>0$,
 $y=0$, and  $f(0)=0$. Let $u_{k}=e^{k\pi i/3}$ and $S=\{t u_{k}:
t\ge 0, k=0,2,4\}$. Let also $0<\d<1$ be such that $\| x\|,\| y\|\le
2$ and $\| x-y\|<\d$  imply that $\|f(x)-f(y)\|<1/2$.

For each $x\in \text{\rm B}_1(0)\cap f^{-1}(0)$ and $k=1,3,5$, use Lemma
4.4 to choose
$\phi_{k,x}:[0,\infty)\to \Rn{2}$ having Lipschitz constant one such
that
$\phi_{k,x}(0)=x$ and $f(\phi_{k,x}(t))=t u_{k}$ for $t\ge 0$.  Let
$D_{k,x}$ be the component of  $\Rn{2}\setminus f^{-1}(S)$ containing
$\phi_{k,x}(0,\infty)$. Noting that
$\text{\rm B}_{\d}(\phi_{k,x}(1))\subset D_{k,x}\cap \text{\rm B}_2(0)$, a
comparison of areas shows that the set of all such $D_{k,x}$ has at
most $N\le 4\d^{-2}$ elements. Suppose now that  $\text{\rm B}_1(0)\cap
f^{-1}(0)$ has infinitely many elements, hence it contains elements
$x\ne y$ such that
$\{D_{1,x},D_{3,x},D_{5,x}\}=\{D_{1,y},D_{3,y},D_{5,y}\}$. Then
$D_{k,x}=D_{k,y}$ for $k=1,3,5$, since the (connected) image of
$D_{k}:=D_{k,x}$ contains $u_{k}$ and so can contain no other $u_{j}$,
and we infer that there are simple curves $\psi_{k}:[0,1]\to \Rn{2}$
such that
$\psi_{k}(0)=x$, $\psi_{k}(1)=y$ and $\psi_{k}(t)\in D_{k}$ for
$0<t<1$. For each pair $k,l=1,3,5$ of different indices, let
$G_{k,l}$ be the interior of the Jordan curve
$(\psi_{k}-\psi_{l})$ (difference in the sense of oriented curves).
If $j\ne k,l$, we note that $G_{k,l}\cap D_{j}=\emptyset$ since
otherwise $D_{j}$ would be bounded. In particular,
$G_{1,3}\cap \partial G_{3,5}=\emptyset$, so either $G_{1,3}\subset
G_{3,5}$ or $G_{1,3}\cap G_{3,5}=\emptyset$. In the former case we
would get a contradiction from $\psi_{1}(0,1)\subset G_{3,5}$, since
$\psi_{1}(0,1)\subset D_{1}$ and in the latter case we would infer
from $\partial (G_{1,5})=\partial (G_{1,3}\cup G_{3,5})$ that
$G_{1,5}\supset G_{1,3}$ and get a contradiction from
$\psi_{3}(0,1)\subset G_{1,5}$.
\hfill\qed

\enddocument